\documentclass[12pt]{article}
\usepackage{graphicx} 

\usepackage{fullpage}
\usepackage{hyperref}
\usepackage{amsmath,amssymb}
\usepackage{orcidlink}

\newcommand{\R}{\mathbb{R}}
\newcommand{\C}{\mathbb{C}}
\newcommand{\bI}{\mathbb{I}}

\newcommand{\tip}{\mathrm{tip}}
\newcommand{\EWS}{\mathrm{EWS}}
\newcommand{\AR}{\mathrm{AR}}
\newcommand{\VAR}{\mathrm{VAR}}

\title{From trends to optimized forecasts: Quantifying the skill and advancing the utility of dynamics-based early warnings for tipping events}

\author{Franco Du Plessis \thanks{University of Exeter, Exeter EX4 4QF, UK} \orcidlink{0000-0002-3064-5550}  \and Victoria Volodina${}^*$ \orcidlink{0000-0003-4490-8777} \and Chris A. Boulton${}^*$ \orcidlink{0000-0001-7836-9391} \and Muhammed Fadera${}^*$ \orcidlink{0000-0002-9882-5005} \and Sneha Kachhara${}^*$ \orcidlink{0000-0002-0418-376X} \and Timothy M. Lenton${}^*$ \orcidlink{0000-0002-6725-7498} \and Paul D. L. Ritchie${}^*$ \orcidlink{0000-0002-7649-2991}\and Peter Ashwin${}^*$ \orcidlink{0000-0001-7330-4951}}

\date{\today}

\begin{document}

\maketitle

\begin{abstract}
In this paper we explore and extend the state of the art for understanding the skill and utility of early warnings of critical transitions (tipping points) in forced nonlinear systems. We highlight some of the challenges and opportunities of transforming estimates of a dynamics-based early warning system, based on trends in stability or resilience into forecasts with high skill. We highlight the importance of (a) quantifying consequences of possible actions in response to warnings that may be false negatives or positives (b) considering finite time horizon predictions to give verifiable predictions (c) assumptions necessary for valid extrapolations of trends. We evaluate approaches to improving the skill and utility of the forecast.
\end{abstract}

\newpage
\tableofcontents

\section{Introduction}
\label{sec:intro}

There is compelling evidence that increasing global mean temperatures or greenhouse gas cumulative emission levels can cause certain components or `elements' of the climate to cross tipping points \cite{armstrong2022exceeding}. Crossing these thresholds can result in rapid changes to climate states via self-reinforcing feedback to states that are new, potentially very challenging and hard or impossible to reverse. Crossing these thresholds may affect regional land ice, ocean ecosystems, hydrology, vegetation and ocean flows with consequences for society that may prove hard to adapt to \cite{richardson2023}.

For this reason, early warning has potentially high societal value, where there will be a link between the lead time of any prediction and the options available for mitigation or adaptation. 
This can add value to the science of climate projections epitomised by the Intergovernmental Panel on Climate Change (IPCC) \cite{IPCC_AR6}. Only in comparatively recent IPCC assessment reports have the underpinning CMIP (Coupled Model Intercomparison Project) models reached the level of maturity that the necessary physics has been included to model some, but not all, of the known possible tipping elements \cite{tipsci2026}. 

In principle, it may be possible to find 
specific thresholds of emissions or temperature that, when crossed, will lead to self-reinforcing feedbacks and tipping events (such as how knowledge of a river level is a predictor of flooding). However, is not practically possible to use such a simple approach to understand tipping triggered by changes in greenhouse gas emissions or global mean surface temperature. The most prevalent approaches are based on measuring a {\em dynamics-based indicator}, such as {\em critical slowing down} (CSD), from carefully selected time series that contain key dynamic information about the tipping element \cite{Lenton08,Thompson2010FoldTippingPredict,scheffer2012anticipating,kuehn2013mathematical,dakos2024tipping}, and this is still an active area of research \cite{huntingford2026}. 

Tipping can be thought of as an abrupt change in a system's output in response to slow changes in input (or forcing). This means we need to identify the system of interest and what is viewed as an input. Then we are in a position to view this as a response of a non-autonomous dynamical system - also called an open system \cite{ashwin2012tipping}. For the case where the forcing changes slowly enough relative to system timescales, and the noise is small enough, any rapid transitions can be associated with slow passage through bifurcation (bifurcation-induced or B-tipping). Even without a detailed physics-based model, it may be possible to find out something about how close a system is to a B-tipping point through exploiting generic relationships between the forcing and system variability. On approaching such a tipping threshold, the system will recover more slowly from perturbations as the restoring forces or damping feedbacks in the system get weaker. The simplest and most used indicator of stability is based on estimation of the lag-one autoregression ($\AR(1)$) coefficient, though there are other ways of quantifying the `resilience' of a stable state, for example based on distance to basin boundary - see for example \cite{krakovska2024resilience,morr2026computing,smith2026estimating}. 

If noise or internal variability is large enough, we can get transitions even without any bifurcations present (noise-induced or N-tipping). If noise is small but forcing is relatively large, the dynamical state may be some way from fixed-forcing equilibria and this can lead to transitions that are not associated with bifurcations (rate-induced or R-tipping) or lack of a transition even if a bifurcation is crossed (safe overshoot scenarios).
For such cases, we might need additional information such as the anticipated future forcing, the current system state, and the rules of system evolution (equations governing system dynamics): see \cite{chapman2025quantifying,ritchie2026evaluating}. 

A warning system of some impending tipping event will never be perfect - there will be a variety of sources of error in estimating the dynamics that may lead to tipping. This includes measurement errors in sensing, errors in state and parameter estimation and model discrepancy (structural uncertainty), which accounts for the mismatch between the model and the physical process of interest \cite{kennedy2001bayesian, volodina2021importance}. Moreover, we argue that it is not possible to quantify the skill of a prediction with no timescale, and so we are concerned with a prediction that we will cross a tipping point within some fixed horizon (i.e. lead) time.

For robust decisions that take such warnings into account, we need to be mindful of the consequences not just of a prediction missing a tipping event (a false negative), but also of a predicted event not coming to pass (a false positive). These consequences will depend critically on the warning or lead time that we are able to give for any prediction. Such approaches are used, for example, to evaluate the prediction accuracy of a receiver at the end of a noisy communication channel and to predict extreme weather events. 
Hence, the skill of a warning system needs to be evaluated in relation to the likely cost/benefit of applying the warning systems, and we are interested in optimising the benefit minus the cost. For example, if some action delivers benefits from true positives and otherwise has no costs, then no early warning is needed to see the benefit of the action. This is, however, rarely the case!

\begin{figure}[ht!]
    \centering
    \includegraphics[width=1.0\linewidth]{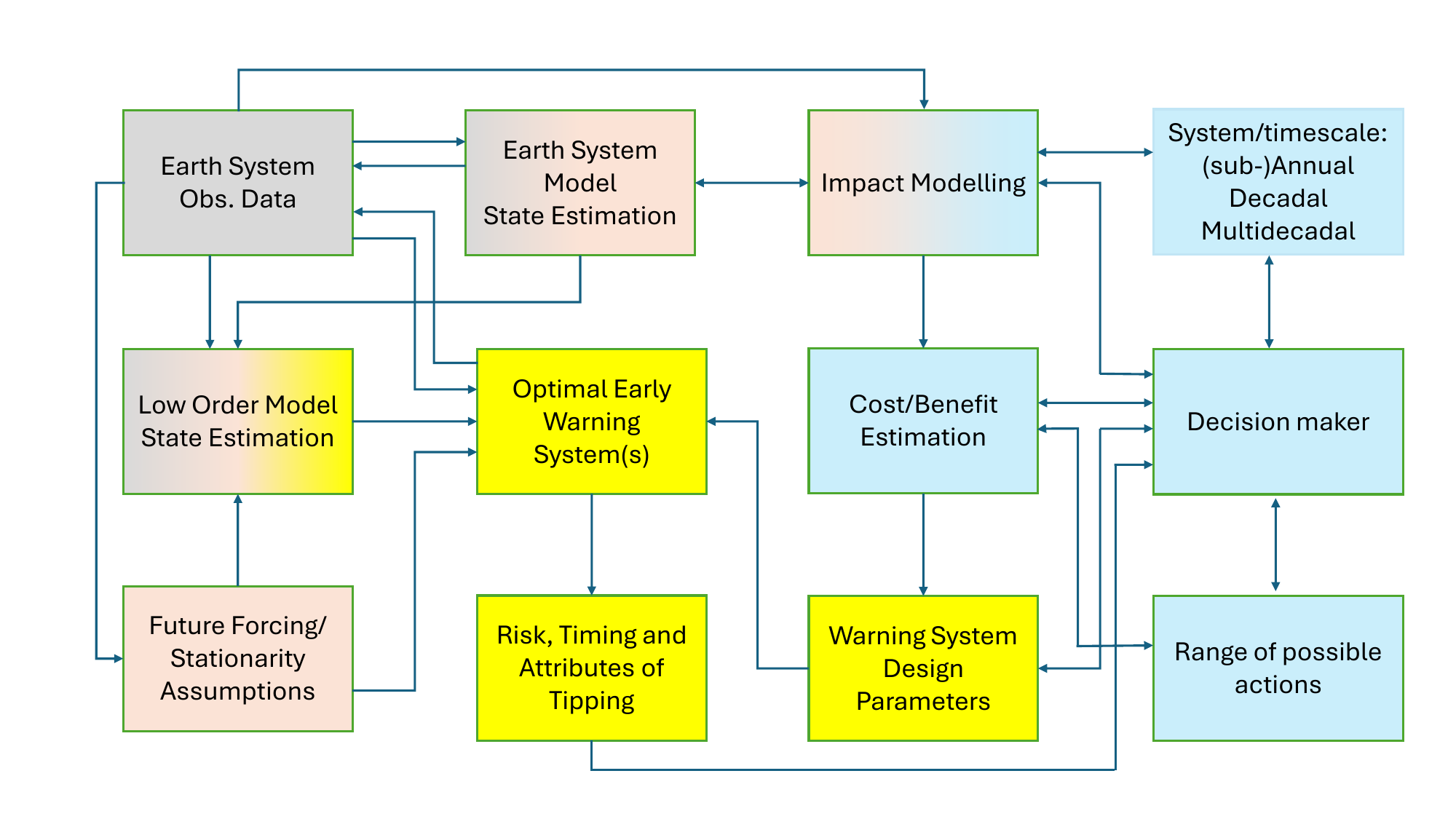}
    \caption{Design of an optimised and skilful dynamics-based forecasting system to warn of approaching Earth system tipping events. A good understanding of potential impacts of events and of predictions of events (and errors in predictions) is important in deciding what should be predictable and what additional information can improve predictive skill. The workflow combines outputs from observations (grey) and models (pink) to produce and optimise an early warning system (yellow) adapted to a specific decision-maker (blue). Note that each decision-maker may have different design parameters that lead to warnings at different thresholds.}
    \label{fig:skill-blobs}
\end{figure}

Figure~\ref{fig:skill-blobs} highlights key aspects that need to be taken into account when designing an optimal early warning system (centre) that optimises the skill of prediction of a climate tipping point. In addition to Earth system measurement and modelling aspects, a simplified low-order model (left) is needed to estimate how close we might be to a tipping point and whether we are likely to arrive there in a given time horizon. To evaluate the skill, we need to understand the implications of possible decisions on impacts (including ``do nothing") to make assumptions about future forcing, system stationarity and impacts of predictions. Note that for multiple stakeholders with varying lead times and action and timescales, we will typically need a ``dashboard" tuned to be optimal for a range of different design parameters. 

We proceed to give two simple illustrative examples. If an action is available with zero cost and only benefit if there is a tipping, then any sort of warning of tipping is superfluous - one should simply implement that action. On the other hand, if there is a large cost of false positive warnings (for example, expensive mitigation actions that are not needed, leading to loss of credibility of a warning system) then we need to weigh this cost against the benefit of actions in the case there is a tipping.

The purpose of this paper is to give an overview of approaches to quantifying and understanding the skill of early warnings for tipping events, and to identify possible steps to improve and optimise such early warnings. We highlight some of the mathematical and decision-making challenges for the remainder of this section. Section~\ref{sec:estimating-trends} looks at approaches to estimating stability and trends in stability, while Section~\ref{sec:uq-forecasts} looks at approaches to uncertainty quantification and skill estimation for dynamics-based early warnings of tipping.
This is followed in Section~\ref{sec:improving-skill} with a look at approaches to improving and optimizing the skill of dynamics-based early warnings in the future, both through better use of available data, and through targeted acquisition of new data.

\subsection{The mathematical challenges}
\label{sec:math-challenges}

The main feature of a dynamics-based early warning is that it tries to estimate the stability of a dynamical state of a (stochastic) nonlinear system, so that trends in the stability can be measured and extrapolated to predict loss of stability at an approaching bifurcation.
While successful in several contexts, standard CSD based methods need several hypotheses to hold, and these are not always made explicit. Violating one or more hypothesis can result in false negatives \cite{Boerlijst2013SilentCollapse} and false positives \cite{rietkerk2025ambiguity,lehnertz2024time}. Moreover, it is often not clear how to check whether the system considered satisfies the hypotheses. 

More precisely, CSD estimations of trends in stability assume the following hypotheses which impose conditions on the predictive ability when tipping is bifurcation-induced:
\begin{itemize}
    \item[\bf H1] Tipping is associated with slow passage through a bifurcation of a subsystem, where a stable equilibrium loses stability
    \item[\bf H2] Noise is relatively simple (typically white/Wiener) and low enough amplitude that linear approximation is valid away from tipping
    \item[\bf H3] There is sufficient data of high enough quality, sampled at high enough frequency to estimate eigenvalues
    \item[\bf H4] The system is in a region where a bifurcation normal form approximation is valid, and the forcing changes not only slowly but with consistent trend
\end{itemize}
A clear challenge is to find early warning predictors that work for a wide variety of tipping elements, and this will be through advancing the mathematical theory to ensure we can improve skill and robustness \cite{huntingford2026}. This is a challenge both in a probabilistic and a deterministic sense. Unpacking these requirements, we list a number of hypotheses we may wish for.

With regard to {\bf H1}, the typical case considered is for a saddle-node or fold bifurcation of an equilibrium, but (a) other bifurcations (in particular Hopf bifurcations) are just as likely to be relevant and (b) the base state may be not equilibrium - it may be periodic or have long-timescale natural variability. Moreover, there may be a significant delay between crossing a threshold and a significant effect arriving. 
With regard to {\bf H2}, the noise may be more structured and much less easy to separate from the nonlinear dynamics, and it may be that long-timescale chaotic variability is a better model for e.g. ocean variability.
With regard to {\bf H3}, estimating stability, together with uncertainty of an equilibrium state from fluctuations can be done relatively simply for a one-dimensional system \cite{Ditlevsen2023AMOCPrediction}, but there are much higher requirements in terms of data for estimating the spectrum of eigenvalues in higher dimensional systems. In particular, any estimates of changes in stability need some uncertainty quantification to suggest that any trend is significant.  We examine this in Section~\ref{sec:estimating-trends}.
With regard to {\bf H4}, if we wish to predict the point in time where a tipping occurs, we need to make assumptions about the validity of a low-order truncated model for the behaviour and its extrapolation to the point where stability is lost \cite{ashwin2025early}. 

\subsection{The decision-making challenges}

In the paper {\em ``Do early warning signals of tipping points lead to better decisions''} \cite{diekert2025doearly}, the authors demonstrate that designing an Early Warning System needs more than just a skilful warning signal - it requires a decision-maker who can take action according to the signal, either by attempting adaptation to minimize damages, or attempting a mitigation to intervene in the system so as to prevent some damages occurring. 

From a decision-maker's perspective this requires a portfolio of prospective planning and robust adaptation portfolios spanning a diverse set of possible solutions \cite{biesbroek_adaptation_2025}, flexible enough to allow responses to change as conditions evolve, since decisions taken today may foreclose options that would otherwise be available in the future. Biesbroek et al. \cite{biesbroek_adaptation_2025} further argue that monitoring the system at risk of tipping is only one part of the problem, since the broader societal context, including shifts in political will, financial flows, and institutional readiness, must also be tracked, as the window in which adaptation or mitigation is viable depends on both the natural system and the human system responding to it. Another complication is that intervention itself can feed back into the EWS. Dai et al.\ \cite{Dai2015StabilityvsResilience} show that systems subject to multiple forcings acting in opposing directions can exhibit EWS failure, and by extension, mitigation policies that simultaneously increase stability by one metric while reducing the basin of attraction by another could place the system in a regime where standard indicators are uninformative or misleading. Barrett and Dannenberg \cite{barrett2022sensitivity} and Hurlstone et al. \cite{hurlstone2025threshold} highlight evidence that increased uncertainty in thresholds is likely to reduce the level of any cooperative action.

Which actions are possible will depend on the timescale, utility and skill of a warning, and so the early warning system design parameters illustrated in Figure~\ref{fig:skill-blobs} may be quite complex. One can make a number of hypotheses for options facing the decision-maker to be evaluated:
\begin{itemize}
    \item[\bf H1D] The set of future possible scenarios is known (to some degree).
    \item[\bf H2D] The range of possible actions one could take is known.
    \item[\bf H3D] The costs/benefits of taking/not taking action in response to an Early Warning can be quantified with some confidence.
\end{itemize}
These are much more likely to be known for the near future, so increased caution is inevitably needed for warnings over longer timescales.

With regard to {\bf H1D}, future scenarios may become increasingly unreliable after tipping occurs as the system may enter states where models have not been subjected to rigorous tuning. In this sense, palaeoclimate models and data are particularly valuable to give analogies, for example, of periods where land ice sheets and ocean circulation patterns have been substantially different.
With regard to {\bf H2D}, one can assume there is a range of known mitigation and adaptation options, but other options may appear based on better knowledge of the system, the impacts, and possible control mechanisms.
With regard to {\bf H3D}, the costs and benefits of taking/not taking action will often be different to different populations and sectors, and so this will tend to be sensitive to the motivations of the decision-maker.

\subsection{Modelling via a low order nonautonomous stochastic systems}

Systems can be mathematically modelled from either physical principles or from observation of its behaviour, specifying how its state evolves in time. An alternative viewpoint is to specify how functions of the state evolve in time. 
As a detailed modelling of the (relatively) rapid processes leading to variability (for example, atmospheric and ocean turbulence in climate models) one can take the approach of Hasselmann \cite{Hasselmann1993optimal} that these are represented by stochastic terms. Suppose the system of interest consists of a coupled system of macroscopic (vector valued) variables $\mathbf{x}$ and microscope variables $\mathbf{z}$
\begin{equation}
\label{eq:FSODE}
\begin{aligned}
\frac{d}{dt}\mathbf{x}&=f(\mathbf{x},\mathbf{z},t)\\
\epsilon \frac{d}{dt}\mathbf{z}&=g(\mathbf{x},\mathbf{z},t)
\end{aligned}
\end{equation}
where $\epsilon$ is a timescale separation. 

If $\epsilon$ is small so that the microscopic variables evolve rapidly and chaotically, one can approximate (\ref{eq:FSODE}) 
with a stochastic differential equation model (also called a Langevin equation) with state $\mathbf{x}(t)$ at time $t$, we write
\begin{equation}
\label{eq:SDE}
d\mathbf{x}(t)  = f(\mathbf{x}(t),t) dt + \sigma(\mathbf{x}(t),t) d\mathbf{W}_t
\end{equation}
where $d\mathbf{W}_t$ is a vector of independent standard white noise processes, $f$ is the ``deterministic" part of the evolution and $\sigma$ is a ``noise strength". The functions $f$ and $\sigma$ are usually assumed smooth in $\mathbf{x}$, and slowly varying with $t$, and $\sigma$ is a matrix such that $\sigma^T\sigma$ is positive semidefinite.

Understanding if and when (\ref{eq:SDE}) is a good approximation of (\ref{eq:FSODE}) are challenging questions that affect early warning skill when it breaks down: we mention some literature on deterministic homogenisation which discusses limiting cases where the reduction may or may not be valid \cite{kelly2017deterministic,melbourne2011note,gottwald2017stochastic}.
Note that the Mori-Zwanzig operator \cite{Zwanzig1961Memory,Chorin2000MZMemory} explicitly splits dynamics into resolved and unresolved components, giving a generalised Langevin equation with explicit memory term (non-Markovian). This can be seen as a way to get a more accurate low-dimension approximation, though not in the form (\ref{eq:SDE}).  
Note that even if (\ref{eq:FSODE}) is known, the state $\mathbf{x}(t)$ is typically only partially observable: we rely on some observation $\mathbf{y}(t)=\phi(\mathbf{x}(t))$ from the system where the observation $\mathbf{y}$ is a vector of typically much lower dimension than the state $\mathbf{x}$. 


\section{Estimating stability and loss of stability}
\label{sec:estimating-trends}

\subsection{Autoregression (critical slowing down) based methods}

The most widely used approach, monitoring what is often called critical slowing down (or CSD) exploits the fact that, under suitable conditions (outlined in Section~\ref{sec:math-challenges}) the stochastic dynamics near a stable equilibrium reduces to a form where stability properties can be estimated directly from time series using autoregressive models.

To make this precise, consider equation \eqref{eq:SDE} for $\mathbf{x} \in\R^N$ with constant additive noise and an explicit time-varying control parameter $p(t)$
\begin{equation}
    d \mathbf{x} = \mathbf{f}(\mathbf{x}, p(t)) dt + \sigma d\mathbf{W}_t,
    \label{eq:SDEnonaut}
\end{equation}
where $\sigma >0$ is a constant noise amplitude and $\mathbf{W}_t$ is an $N$-dimensional standard Wiener process \cite{gardiner2004handbook}. If $p(t)$ is slowly varying then we can approximate (\ref{eq:SDEnonaut}) as the autonomous, quasistatic or frozen system
\begin{equation}
    d \mathbf{x} = \mathbf{f}(\mathbf{x}, p) dt + \sigma d\mathbf{W}_t.
\end{equation}
Assuming {\bf H1} the system possesses a stable equilibrium $\mathbf{x}^*(p)$ satisfying  $\mathbf{f}(\mathbf{x}^*, p) = \mathbf{0}$, and by {\bf H2} the noise amplitude is small enough that the linearisation about $\mathbf{x}^*(p)$ is valid.  Treating $p$ as fixed at the point of linearisation and writing $\mathbf{y} = \mathbf{x} - \mathbf{x}^*(p)$ for small deviations from equilibrium, the dynamics reduce to a linear SDE (an Ornstein-Uhlenbeck equation):
\begin{equation}
\label{eq:linearised-sde}
    d\mathbf{y} = J(p)\mathbf{y} dt + \sigma d\mathbf{W}_t,
\end{equation}
where $J(p) = \partial \mathbf{f}/\partial \mathbf{x}\big|_{\mathbf{x}^*(p)}$ is the Jacobian, whose eigenvalues $\lambda_i$ all have strictly negative real parts for a linearly stable equilibrium. As $p$ approaches a bifurcation  point $p_c$, {\bf H1} and {\bf H4} together imply that the real part of the leading eigenvalue satisfies $\mathrm{Re}(\lambda_1(p)) \to 0^{-}$. 
The goal of autoregressive methods is then to estimate eigenvalue structure and its trend from observed time series, subject to the data requirements of {\bf H3}.

If there is a single real eigenvalue that becomes unstable then, in principle, one can project onto the leading eigenvector to reduce the dynamics to the scalar equation $y\in\R$ for the critical mode
\begin{equation}
dy = -\alpha(p)y dt +\sigma dW_t,
\end{equation}
where $\alpha(p) = -\lambda_1(p) >0$  is the {\em restoring rate}, so that $\alpha(p) \to 0^{+}$ as $p \to p_c$. Discrete time sampling of this process with time interval $t_{n+1}-t_n=\Delta$ gives
\begin{equation}
y_{n+1}=\phi y_{n}+ e_n
\end{equation}
with $e_n$ i.i.d. Gaussian and zero mean, where $\phi=e^{-\Delta \alpha}$. This is exactly a lag-1 autoregressive process $\AR(1)$ with coefficient $\phi = e^{-\lambda \Delta t}$. Hence, fitting an $\AR(1)$ process to observations directly can give an estimate of $\alpha$ via $\hat{\alpha} = -\frac{\log \hat{\phi}}{\Delta t}$, with CSD manifesting as $\hat{\phi} \to 1$ \cite{dakos2008slowing}. However, reality is usually not so simple in that the optimal observable is rarely observable, and even if it is partially observable, there will often be internal noise ``interference'' that can confound the CSD signal \cite{Morr2024InternalNoiseInterference}.

Note that scalar projection does not account for Hopf bifurcations where a complex conjugate pair of eigenvalues approach the imaginary axis simultaneously, producing oscillatory behaviour that a scalar $\AR(1)$ model cannot represent. A more general strategy for \eqref{eq:linearised-sde} works with a leading eigenspace which may or may not be two dimensional, by estimating (at least part) of the spectrum of eigenvalues in terms of a
$\VAR(1)$ process 
\begin{equation}
\mathbf{y}_{n+1} = \Phi\,\mathbf{y}_n + \boldsymbol{\varepsilon}_n
\end{equation}
with $\Phi = e^{J\Delta t}$ and $\varepsilon_n$ a Gaussian noise, and the full eigenvalue spectrum can be estimated  from multivariate time series \cite{Williamson2015MultiDBifDetection,
hobden2026early,lutkepohl2013introduction}, though inevitably estimating more variables to the same level of accuracy will require more data. Note that an important consideration in estimating return rates/eigenvalues is the performance (expected bias and variance) of the estimators with small or imperfect samples \cite{Ditlevsen2023AMOCPrediction,pilipovic2025strang}.

Extensions both for univariate and multivariate approaches exist. This includes more general theoretical frameworks where the local dynamics is estimated not just using a linear SDE (\ref{eq:linearised-sde}) but also including nonlinear term, notably by considering noisy ``normal forms'' near bifurcations \cite{Kuznetsov2004AppliedBifurcationTheory,Berglund2006NoiseInducedPhenomena}. Tipping associated with bifurcation crossing can be classified into a small number of ``typical'' cases - for a single parameter in the absence of special structures this is just the fold (saddle node) bifurcation with normal form for $y\in\R$ with $p_c=0$
\begin{equation}
\label{eq:fold-normform}
dy_t= (ay_t^2 + bp(t))dt+\sigma dW_t
\end{equation}
with $ab<0$ for tipping at $p\rightarrow_-0$, 
and subcritical Hopf bifurcation with normal form written for $z\in\C$ as
\begin{equation}
\label{eq:hopf-normform}
dz_t= ((ap(t)+i\omega)z_t+b |z_t|^2z_t) dt+\sigma dW_t.
\end{equation}
with $a>0$ and $b>0$ for tipping at $p\rightarrow_-0$. Other normal forms considered in the literature include transcritical and pitchfork bifurcations \cite{Kuehn2011CSD,Kuehn2013CSDVar} which can be relevant in the presence of invariant subspaces and symmetries. Normal forms can be particularly useful in that we do not need de-trending, and the scaling of the restoring rates close to bifurcations is a property of the bifurcation. For example, for the fold bifurcation (\ref{eq:fold-normform}), the restoring rate for low noise and slow variation of $p<0$ is given by 
\begin{equation}
    \alpha \sim 2\sqrt{ |abp|}
\end{equation}
and so if $p$ is locally linear in $t$ then $\alpha^2$ will go to zero linearly in $t$ as we approach the bifurcation $p\rightarrow_-0$: see \cite{Ditlevsen2023AMOCPrediction,ashwin2025contrasting}. By contrast, for the Hopf bifurcation (\ref{eq:hopf-normform}) the restoring rate is given by
\begin{equation}
    \alpha \sim |ap|
\end{equation}
and so if $p$ is locally linear in $t$ then $\alpha$ will go to zero linearly in $t$ as $p\rightarrow_-0$.

It is often ignored that the transformation of observables into normal form involves a special choice of nonlinear transformations in state and forcing parameter that make the system ``as simple as possible'', and moreover transformation to normal form is only typically valid close to the bifurcation, where the closeness depends on the size of nonlinear terms that have been transformed away \cite{ashwin2025early}.

Other research considers more general noise processes. If the noise is correlated, \cite{Boettner2022CSDCorrNoise} derive a corrected estimator for $\alpha$ by modelling the noise itself as an OU process. 
Alternatively, higher-order $\AR(\tau)$ models provide a more robust estimate of $\alpha$  by capturing residual autocorrelation structure that would otherwise bias the $\AR(1)$ coefficient.

An extensive list of early warning signal methods based on detecting critical slowing down has been curated by  Dakos et al. \cite{dakos2014critical}. Generally, these are calculated on a moving window of a time series of the system in question. Considerations should be taken into account during their use when choosing the window length used to calculate these indicators. A window length should be long enough to observe the dynamics of the system, but short enough that the resulting indicator time series is long enough to observe significant trends. Robustness analysis should also consider multiple detrending techniques to remove long-term trends that affect these signals.

Bayesian state-space approaches extend this theoretical framework by modelling the $\AR(1)$ coefficient $\phi(t)$ as a latent time-varying process and inferring properties of its posterior trajectory, for example, whether $\phi(t)$ has an increasing trend, directly from data \cite{myrvoll-nilsen_bayesian_2025}: see Section~\ref{sec:uq-forecasts}.

\subsection{Other methods of estimating loss of resilience} 

There has been a lot of recent work on operator-theoretic approaches that transform nonlinear stochastic systems into linear operators acting on (approximations of) function spaces \cite{lucarini2024detecting,chekroun2020ruelle,zagli2026bridging}. In particular, the spectral properties of Koopman operators (for observables) or their duals transfer/Perron-Frobenius operators (for phase space) reveal intrinsic features of a system's response to perturbations, including precursors of loss of stability. 

The eigenfunctions determine coordinates in which the dynamics evolve linearly.
Several approaches are available; these include
Dynamic Mode Decomposition (DMD) \cite{kutz2016dynamic}, Extended DMD \cite{williams2015data}, or other mode decomposition methods \cite{schmid2022dynamic}.
Some state-space formulations estimate continuously varying AR coefficients as latent variables via the Kalman filter \cite{ives2012detecting,kwasniok2024semiparametric}. Probabilistic variants fit a time-varying $\AR(1)$ model using Bayesian inference, returning a posterior distribution over the trajectory of $\phi(t)$ \cite{laitinen2021probabilistic}. A probabilistic time-varying variant of $\VAR(1)$ extends similar ideas to the multivariate setting \cite{laitinen2022probabilistic}. Other approaches include estimating nonlinear terms either in general \cite{kwasniok2024semiparametric} or fitting to an assumed normal form \cite{Ditlevsen2023AMOCPrediction}. Other data-driven modelling approaches that may be useful include equation-free methods \cite{kevrekidis2009equation} or symbolic regression methods such as SINDy \cite{brunton2016discovering}.

In all cases, these methods can result in an indication of a trend of stability. To what extent such a trend can be extrapolated to indicate an impending tipping point remains a vital, but not necessarily simple, question to answer \cite{ashwin2025early,morr2026extrapolation}.
One challenge that is common to all methods of estimating resilience is the need to pick observables that are informative and that have a good signal-to-noise ratio. A good physical knowledge of the system in question will help, but in some cases the only skilful observables may be inaccessible in current datasets.

\subsection{Early warnings for edge states of rapidly forced systems}

In many systems of practical interest, we are interested in forcing that may exceed the capacity of the system to recover, either because the system timescale is particularly short, or the forcing is particularly fast: we call this a rapidly or strongly forced system. In such cases, {\bf H1} will fail to hold, and it is not possible to use threshold crossing to determine the eventual outcome after rapid forcing. Worse still, the system may take states that are far from any invariant sets of the unforced system.

If there are reliable process-based models of the system available, then it may be possible to give early prediction of eventual states by looking at where the trajectory is relative to so-called {\em edge states} \cite{chapman2025quantifying,lohmann2025role,ritchie2026evaluating}. These are dynamically invariant sets that are on basin boundaries between different attractors. These edge states may have distinct, detectable properties that provide early warning, and can help understand the influence of current state, future forcing and future variability on the eventual fate.
In such cases a variety of rate-induced effects can appear. This includes rate-induced tipping, where tipping occurs without any bifurcation \cite{ashwin2012tipping,ashwin2017parameter,wieczorek2023rate}, or safe overshooting where crossing a threshold for a sufficiently short time does not result in tipping.

If a system is rapidly but transiently forced (i.e. such that the forcing is asymptotically stationary) one can construct a time-dependent threshold in phase space (called the R-tipping threshold in \cite{wieczorek2023rate}) that determines the eventual fate of the system. This suggests a practical way to give early warnings for strongly forced systems - to approximate this moving threshold \cite{chapman2025quantifying} and measure the distance from this as a geometric early warning. This can give skilful predictions of final state, as discussed in \cite{ritchie2026evaluating}. In particular, in general, early warnings of tipping for cases of rapid forcing will depend on knowledge of future forcing, not just current state and trend.
We discuss further simulation-based approaches in Section~\ref{subsec:sim_approaches}.

\subsection{Estimating trends in stability}

Many CSD estimators are estimated sequentially over time, either by computing them within a rolling window or by inferring a full trajectory from the estimated quantity. In either case, determining 
whether any observed trend towards instability is genuine rather than arising by chance requires care. Kendall's $\tau$ provides a natural summary statistic for the degree to which an indicator monotonically increases over the observed window, but a statistically significant trend can arise for reasons entirely unrelated to CSD \cite{boettiger2012early, Hastings2010NoEWS} (see also \cite{rietkerk2025ambiguity}). Surrogate methods address this by generating artificial time series that are consistent with a specified null hypothesis while preserving selected properties of the data; this requires both an appropriate null model and a robust test statistic against which the observed behaviour can be evaluated \cite{lehnertz2024time}. The most common approach is Fourier phase randomisation, which preserves  the amplitude spectrum while randomising phase information and thereby maintains the linear autocorrelation structure of the series \cite{Boettner2022CSDCorrNoise, dakos2008slowing}. Block bootstrapping offers an alternative that instead preserves the temporal dependence structure by resampling contiguous blocks of observations \cite{dakos2008slowing}. The choice of surrogate model is however non-trivial. Phase-randomised surrogates may not always constitute an appropriate null for CSD, particularly in the presence of data gaps or when the underlying dynamics are better represented by low-order autoregressive processes, where Ben-Yami et al. \cite{BenYami2023AMOCDataCSD} argue that an AR(2)-based surrogate might be superior.

\subsection{Forecasting from deep learning}
\label{sec:deep-learning}

Classification-based deep learning (DL) models have recently been used to provide early warning of bifurcation-induced  \cite{bury2021deep, deb2022machine}, rate-induced \cite{huang2024anticipating} and noise-induced tipping \cite{zhang2025anticipating}. DL models are trained on large datasets derived from normal forms of bifurcations and from mathematical models exhibiting rate-induced tipping. The flexibility of DL models enables them to learn higher-order statistics that CSD-based methods may otherwise overlook. This allows them to learn tipping from non-tipping with a higher true positive rate and a lower false positive rate.

Although DL models are generally considered black boxes, methods from explainable AI, such as attention scores \cite{boulton_early_2026} or layer-wise propagation \cite{huang2024anticipating, zhang2015predictability}, can be used to understand, for example, when the model begins to detect early warning of tipping. While the current models are very good at providing early warning with potentially longer lead times than critical slowing down \cite{bury2021deep}, sometimes predicting the probability of tipping equal to $1$ several lead times before tipping, they remain prone to false positives, especially for cases outside their training set \cite{dablander2022deep}. Moreover, the interpretation of the probability of tipping being $1$ before tipping in the context of varying future forcing is unclear. Thus, alternative early warnings which are not classification-based, such as distance to tipping in the forcing \cite{guardamagna2026predicting} or time-to-tip \cite{boulton_early_2026}, are being explored. It is nevertheless crucial to have a representative set of nulls for training these classification-based models for use in an operational early warning system, where the cost of acting on false positives could be high. 
This falls into the broader issue of generalisation for deep learning models. They can be very useful in generalising to unseen data similar to their training set, but struggle on the same task outside their training data. This has prompted the use of in-context learning \cite{sevinchan2026context} based on Bayesian inference, which allows deep learning models trained on bifurcation-induced tipping to generalise to rate- and noise-induced tipping with reasonable skill.  

DL models should also be thoroughly benchmarked across different preprocessing procedures, including those performed before CSD-based analysis. In particular, there is an argument to be made that leaving the trend can help these models develop better skills in contexts where the system is far from a quasistatic equilibrium \cite{huang2024anticipating, zhang2025anticipating}. However, detrending could also be considered a form of normalisation, a standard procedure in most deep learning training routines \cite{batchnorm, santurkar2018does}.

\subsection{Simulation-based approaches}
\label{subsec:sim_approaches}

There has been some notable work in simulating tipping events in climate models ranging from an idealized box model proposed by Stommel \cite{stommel61} to a state-of-the-art global climate model (GCM)\cite{van2023asymmetry}. In these experiments, the AMOC tipping and/or weakening is induced by an external forcing such as increases in greenhouse gases \cite{bellomo2021future} or freshwater input \cite{van2023asymmetry}. Due to computational complexity of climate models, these experiments have been performed using single-member runs with the aim to define the threshold on forcing beyond which the current AMOC state may become unstable (``safe operating-space'') \cite{van2023asymmetry} as well as analyse the climatic response to the slowdown. \cite{mehling2024limits} and \cite{cini2026assessing} highlight some of the limitations of the single-member runs approach such as long transient behaviour, ghost states, rate-induced effects as well as internal variability, none of which can be diagnosed from a single realisation. In addition, these models suffer from various uncertainties mentioned in Section \ref{sec:intro}, that need to be carefully quantified and synthesised when we relate these models to the physical system of interest. \cite{caiado2015bayesian} adopt approaches from Bayesian statistics, emulation and history matching, to characterise which parameter configurations lead to tipping in the four-box model of the Atlantic. Accounting for observation error, emulator uncertainty and structural discrepancy, the analysis returns not a single critical parameter value but a region of parameter space corresponding to AMOC collapse.

Recently, rare event simulations, which are algorithms specifically designed to sample low-probability events, have been adopted as an alternative to deterministic predictions. In particular, \cite{castellana2019transition} estimated the probability of an AMOC transition due to noise in freshwater flux in a conceptual box-model. In this study, two types of transition are considered, namely partial and full collapse of the AMOC, with the probability of a partial collapse over the next 100 years estimated at 0.15. \cite{jacques2024estimation} perform similar analysis by combining a rare-event algorithm, trajectory-adaptive multilevel splitting (TAMS), in which trajectories are ranked by a score function obtained from a machine-learned approximation of the committor. This modification allowed to increase the efficiency and accuracy of the rare-event algorithm. \cite{cini2026assessing} also use TAMS to estimate the probability of a strong AMOC weakening under different CO\textsubscript{2} concentration levels and Shared Socioeconomic Pathway (SSP) 
scenarios in PlaSim-LSG, a general circulation model of intermediate complexity. In addition to providing probability assessment of tipping events, rare event algorithms have been successfully adopted in simulating AMOC noise-induced tipping driven solely by internal climate variability \cite{cini2024simulating}.

\section{Uncertainty quantification and skill estimation for forecasts}
\label{sec:uq-forecasts}

We are concerned with when can a trend in stability actually be used as a prediction or forecast of a future event. For a forecast to be testable, we need to be able to evaluate whether or not a prediction is true or false. A prediction that tipping will or will not occur at some unspecified time in the future is not testable in any real sense. For this reason it is appropriate to give predictions that are over a finite time horizon. Two advantages of this approach are (a) extrapolation of trends will typically only be possible over sufficiently short timescales and (b) it facilitates the evaluation of consequences of possible actions in response to predictions.

\subsection{Extrapolating trends to forecasts}
\label{subsec:extrap-trend}

\cite{ashwin2025early} illustrates the importance of extrapolation for early warnings and some of the issues that can cause extrapolation to fail. Suppose at time $t = t_p$, we wish to predict whether tipping is going to happen for $t\in[t_p,t_p + \Delta]$ for some time horizon $\Delta>0$. We say a \emph{tipping event is underway} at time $t$ if some scalar observable $\Psi$ of the system (typically, the rate of change of some state variable) that leaves a region of linear stability and undergoes an abrupt change passes a chosen threshold $\theta$. 

If we define an indicator function
\begin{equation}
I_{\tip,\Delta}(t_p) := \bI_{\min\{s\geq 0~:~\Psi(x(t_p+s))\geq \theta\} \in [0,\Delta]}
\end{equation}
so that $I_{\tip,\Delta}(t_p)=1$ if there is tipping at some point in the time interval, then an {\em early warning predictor} is a causal function $I_{\EWS,\Delta}(t_p)\in\{0,1\}$ that is a predictor for $I_{\tip,\Delta}(t_p)$. 

One can use a {\em scoring classifier} \cite{fawcett2006introduction} to define such a predictor. This is a real-valued function $C_{\Delta}(t_p)$ such that
\begin{equation}
I_{\EWS,\Delta,\kappa}(t_p) :=\bI_{C_{\Delta}(t_p)\leq \kappa}
\end{equation}
The threshold $\kappa$ for the classifier can be used to balance the ability to correctly detect negatives or positives. As noted in \cite{ashwin2025early} there are limiting thresholds where $I_{\EWS,\Delta,\kappa}$ systematically predicts no tipping or inevitable tipping. 

\begin{figure}
    \centering
    \includegraphics[width=0.9\linewidth]{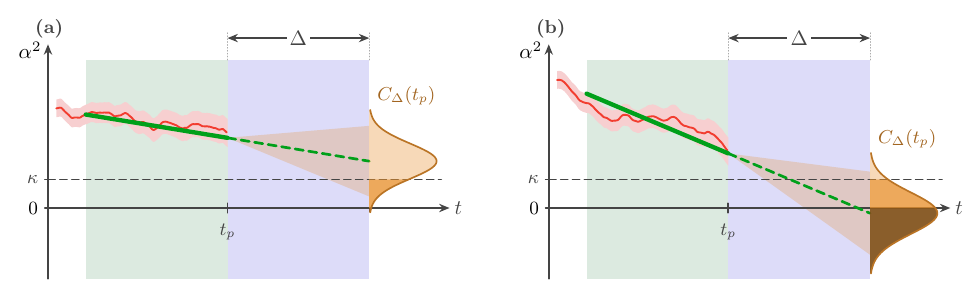}
    \caption{Construction of a  scoring classifier $C_{\Delta}(t_p)$ to predict a tipping point in time horizon $\Delta$ from a prediction time $t_p$, assuming proximity to a fold bifurcation. We assume an estimator (red) for $\alpha^2(t)$ with confidence interval (light red) up to prediction time $t_p$. Linearly extrapolating the green line from the light green region to $t_p+\Delta$, gives a classifier with a distribution that predicts tipping occurs if it is below some threshold.  For threshold $\kappa$ we predict tipping in the blue region $[t_p,t_p+\Delta]$ if $C_{\Delta}(t_p)<\kappa$; By taking $\kappa$ non-zero we can bias towards avoiding false positive or negatives. For more details, see \cite{ashwin2025early}.}
    \label{fig:ews-extrapolation}
\end{figure}

Moreover, given a trend in stability up to time $t_p$, one might hope to construct a scoring classifier that is a good predictor, by extrapolating that trend to predict whether instability occurs within the given time window. Figure~\ref{fig:ews-extrapolation} illustrates the construction of the classifier $C$ via extrapolation.

To illustrate how these definitions can be applied in practice, let's consider the AMOC fingerprint, for instance, sea surface temperature (SST), at time $t_p$, denoted as $x(t_p)$. An observable scalar $\Psi$ here denotes a rate of change in the AMOC fingerprint that exceeds a threshold, $\theta$, within a horizon $\Delta$ (e.g., 20-50 years). We can use the return rate $\hat{\alpha}$ estimated from an observation window as an early warning indicator. \cite{ashwin2025early} demonstrated how to obtain the family of early warning predictors by extrapolating $\hat{\alpha}(t)$. We can then construct a scoring classifier with uncertainty coming from the extrapolation fit:
\begin{equation}
   C_{\Delta}(t_p) := \beta_0(t_p) + \beta_1(t_p)\Delta. 
\end{equation}
The threshold $\kappa$ for the early warning predictor can be chosen to minimise both false negatives and false positives by testing the predictor on the model runs, simulated experiments and/or observational data that contain tipping and non-tipping trajectories. The uncertainty about predictions can be represented using one of the approaches presented in Section \ref{subsec:UQ}.

\begin{figure}
    \centering
    \includegraphics[width=0.9\linewidth]{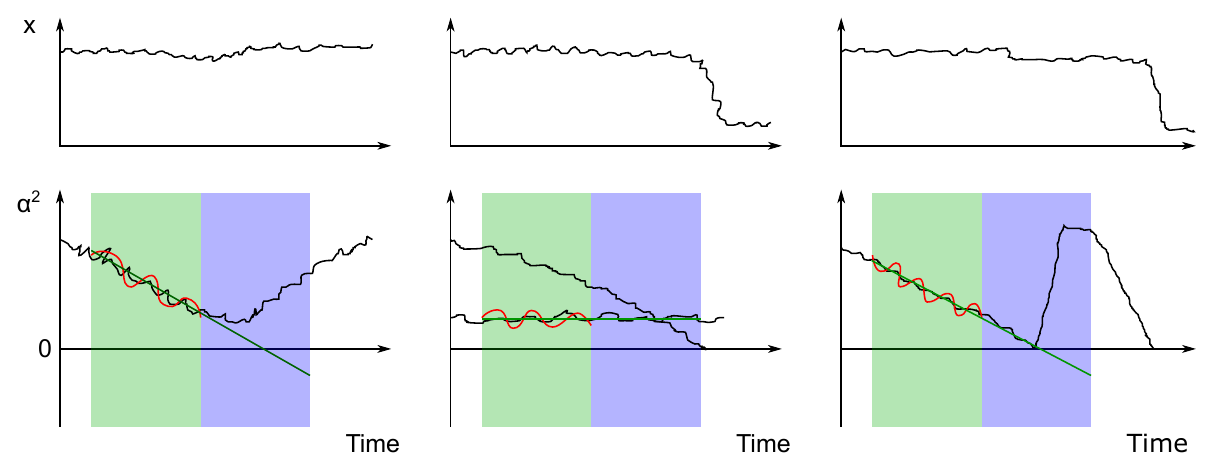}
    \caption{Potential pitfalls with extrapolating trends for early warning from a signal $x$ (top panels). The bottom panels show extrapolation of an estimator for $\alpha^2$ from the end the green region into the blue region. We see that extrapolation might fail because of non-monotonic forcing (left column), a change in leading eigenvalue (middle column), being outside the region of validity of the normal form for the later tipping point associated with a large changes in $x$ (right column).}
    \label{fig:ews-extrapolation-issues}
\end{figure}

As with any extrapolation, this is liable to become increasingly invalid as time proceeds, as discussed in \cite{ashwin2025early}. This paper highlights that the extrapolation can break down for a number of independent reasons, even for cases of B-tipping, see Figure~\ref{fig:ews-extrapolation-issues}:
\begin{itemize}
    \item The forcing may be a nonlinear function of time and  may reverse direction (as in the left column of Figure~\ref{fig:ews-extrapolation-issues}), or have a sudden change of speed (for example, in a system that is downsteam of a system that undergoes tipping).
    \item There may be eigenvalue crossings such that the trend of the leading eigenvalue changes quite abruptly (middle column of the Figure~\ref{fig:ews-extrapolation-issues}).
    \item The extrapolation to a tipping point is based on the assumption that the system is in a region of validity for the normal form. This may not hold unless the system is very close to the tipping, or there could be a much smaller tipping event prior to a larger one that would also fail to be detected (right column in Figure~\ref{fig:ews-extrapolation-issues}). It will also depend on the type of bifurcation being approached.
\end{itemize}

\subsection{Quantitative measures of skill}
\label{sec:skill-measures}

An ROC (receiver-operator characteristic) analysis \cite{fawcett2006introduction} can be used to quantify the rate (i.e. proportion) of false positives $FPR$ and the rate of true positives $TPR$ as a function of the threshold $\kappa$. Note that a variety of other notations can be used to characterise skill; these include 
$$
    \mathrm{specificity}=1-FPR,~~\mathrm{ sensitivity}=TPR.
$$
This approach has been used to understand skill of EWS, for example in:  \cite{Proverbio2022COVIDEWS,Brett2020ROCEWS,zhang2025anticipating}; see \cite{lehnertz2024time} for a comprehensive review with a particular emphasis on surrogate methods.

On varying $\kappa$, note that the ROC goes from $(FPR,TPR)=(0,0)$ to $(1,1)$; for too low $\kappa$ the predictor always predicts no tipping, whereas for too high $\kappa$ the predictor always predicts tipping. An {\em optimal choice} is the value $\kappa=\kappa^{Opt}(t_p,\Delta)$ that corresponds to the closest point to $(0,1)$, i.e. that minimises both false negatives (high sensitivity) and false positives (high specificity). Clearly, this will depend on $t_p$, $\Delta$ and the ensemble considered. 

One particularly useful measure of the optimality is expressed by the {\em area under the curve} (AUC) for this ROC \cite{fawcett2006introduction,lehnertz2024time}; in that if AUC is close to one, we can choose a threshold that simultaneously gives very few false positives and almost all true positives. 
However, it should be noted that AUC is  sensitive to the relative proportions of tipping and non-tipping examples in the evaluation set, and thus can give a false impression under imbalanced classes \cite{lehnertz2024time}. Surrogate-assisted ROC analysis provides one route to addressing this \cite{kuhlmann_seizure_2018, lehnertz2024time}.
We note that other criteria can be used to assess prediction skill, and one can also weigh by cost of false/true negative/positives.

Once an optimal threshold has been selected, the \textit{informedness}, calculated as specificity + sensitivity - 1, can be used to measure the performance of an EWS \cite{ashwin2025early}. A value of one would indicate a perfect categorisation of tipping outcomes, whereas a value of zero would be equivalent to random chance and -1 always provides the wrong categorisation. Informedness therefore offers an improvement over using accuracy, as it is independent of the relative proportions of tipping and non-tipping examples in the evaluation set.    

Beyond discrimination, a fully probabilistic framework enables the use of calibration metrics and proper scoring rules. Calibration metrics such as expected calibration error (ECE) and calibration plots can be used to assess the empirical frequency with which observed values fall within the nominal prediction interval. However, a well calibrated prediction may not lead to a useful uncertainty quantification. To access both calibration and concentration, we can use proper scoring rules \cite{gneiting2007strictly}. Measures such as the Brier score \cite{ferro2012bias} provide a more complete picture of skill than AUC alone. 

Machete \cite{machete_early_2013} develops this idea explicitly in the context of early warnings, showing that mixing conditional density forecasts with the unconditional density can improve both calibration and sharpness, and arguing that this probabilistic framing of skill is more informative than binary detection metrics when the model is subject to misspecification.

\subsection{Uncertainty quantification}
\label{subsec:UQ}

In this section, we review approaches to uncertainty quantification for methods in estimating stability and trends in stability from Section  \ref{sec:estimating-trends}. A decision-maker weighing the cost of acting against the cost of not acting needs a notion of how likely each outcome is. For an early warning system, these likelihoods and the uncertainties underlying them arise from the combined effect of three stages, namely 
\begin{itemize}
\item Construction of an observed time series
\item Estimation of stability from that time series
\item Extrapolation of that estimate across a forecast horizon. 
\end{itemize}
We attempt to address each stage in turn. 

A frequently overlooked source of uncertainty sits upstream of the estimation methods of Section~\ref{sec:estimating-trends}, in the observation-to-data pipeline itself. Observational datasets are the output of processing chains involving aggregation, calibration, interpolation, infilling, smoothing, and resampling, and each of these steps can modify the variance and autocorrelation structure from which CSD indicators are computed \cite{Smith2023ReliabilityResilienceEstimation}. A change in an indicator can therefore reflect the processing as well as the underlying system. This uncertainty can be explored by repeating the analysis across defensible dataset constructions, and the resulting spread in forecasts provides one measure of the sensitivity to observational choices \cite{benmayi2024uncertainties}. Such variation is not the same as sampling uncertainty, but it provides a way of exposing a source of uncertainty that would otherwise remain hidden. The QA4EO framework provides a template for making uncertainty associated with the construction of the observed times series visible and traceable through the observation-to-data pipeline \cite{woolliams_general_2025}.

The sources of uncertainty in an early warning are often dominated by the assumptions underlying CSD-based trend estimation, which are set out in Section \ref{sec:math-challenges}, and these are among the least amenable to direct quantification. If these assumptions fail, the estimated trend may no longer represent the quantity it is taken to measure, and an uncertainty calculated conditional on the assumed model cannot capture the resulting structural error \cite{ashwin2025early}. Since there is generally no accepted way to quantify the probability that these assumptions hold, a forecast is always conditional on them, and this should be stated directly. In addition, a structural discrepancy assessment \cite{goldstein2025systematic,goldstein2013assessing} can be performed, where these assumptions can be treated as sources of internal discrepancy and the effect of these simplifying assumptions made can be investigated by direct experimentation.

Many of the approaches in Section~\ref{sec:estimating-trends} fit statistical models, either locally within rolling windows or to the record as a whole, and parameter uncertainty therefore enters naturally. This is one of the more tractable sources, since estimation theory provides ways of quantifying the sampling variability of fitted parameters (see \cite{gelman_bayesian_2013} and \cite{cox2006principles}). In applications, these statistical models can suffer from misspecification since the assumptions of a well-specified prior and model rarely hold in practice \cite{matsubara2022robust}. Various approaches have been proposed to address this issue from a Bayesian perspective, including generalised Bayesian inference \cite{bissiri2016general}, that adopts flexible loss functions and divergence measures in place of the traditional likelihood functions. The frequentist alternatives include M-estimation \cite{stefanski2002calculus}, and empirical likelihood \cite{owen2001empirical}. Conceptual and climate models are increasingly used to simulate the conditions under which tipping occurs (see Section \ref{subsec:sim_approaches}). In addition to parameter uncertainty and observation error, this type of models also suffers from model discrepancy, which can be accounted for by introducing a model-discrepancy term, treated as an additional uncertain component representing the systematic differences between the model and the observed system \cite{kennedy2001bayesian}. Its magnitude and structure must be constrained by prior knowledge, since a discrepancy term flexible enough to absorb the signal of interest becomes confounded with the model parameters themselves \cite{brynjarsdottir_learning_2014}. This does not remove the risk of misspecification, but it allows the reported uncertainty to reflect the possibility that the assumed model is only an idealised approximation.

Within the frequentist framework \cite{cox2006principles}, uncertainty in the forecast quantity can be expressed as a confidence interval, obtained either from asymptotic theory or by resampling procedures such as parametric bootstrapping. Ditlevsen and Ditlevsen \cite{Ditlevsen2023AMOCPrediction}, for example, simulated trajectories from their fitted model and repeated the estimation procedure to obtain a bootstrap confidence interval for the tipping time. In a Bayesian analysis the posterior (or posterior predictive) distribution provides the corresponding probabilistic description, so the probability of tipping within a specified forecast horizon is the posterior mass assigned to that event (see \cite{myrvoll-nilsen_bayesian_2025} as an example). When dealing with conceptual and climate models, we can perform forward uncertainty quantification and propagate major sources of uncertainty through the model to investigate their effect on the model outputs of interest \cite{kimpton2025challenges}. Forward uncertainty quantification replaces a single deterministic forecast with a distribution for the forecast quantity, from which statistical summaries, such as mean, median and variance, as well as the probability of any particular outcome can be computed.

The reliability of such uncertainty estimates cannot be established for real climate tipping events, since the outcome is known only once the transition has occurred. Testing therefore relies heavily on model output and designed ensembles in which both tipping and non-tipping trajectories are represented; these ensembles may be from model runs, observational data sets or, for example, simulations of observations via an observing system simulation experiment (OSSE). Both tipping and non-tipping examples are required, since evaluating an EWS only on systems that tipped measures how often a warning precedes a transition but says little about how often the same warning would have been produced when no transition follows. Conditioning on observed transitions in this way can produce the prosecutor's fallacy and inflate the apparent skill of a warning system \cite{boettiger2012early}. The same problem is considered by \cite{rietkerk2025ambiguity}, who instead propose to adopt Bayes' theorem with a prior probability of a tipping point approaching to obtain the required probability. This prior can be derived from several sources, including model simulations and palaeoclimate data. In addition, expert elicitation can be adopted to gather knowledge and beliefs of subject-matter experts about tipping points. 
\cite{lam2022climate} focuses on elicitation protocols proposed by \cite{Lenton08} and \cite{kriegler2009imprecise} and highlights that expert elicitation here does more than quantify uncertainty; it also serves to detect and characterise tipping points and their interactions. 

 The role of uncertainty quantification is therefore not to produce a single number that removes the ambiguity of tipping, but to make explicit how the available evidence, modelling assumptions and data limitations combine to support a decision.

\subsection{Tipping duration, timing and decision impacts}
\label{subsec:tipping-duration-timing}

A skilful prediction, at time $t_p$, that tipping will start within some time horizon $\Delta$ is only part of the story; we would like to know the likely timing of $t_{\tip}$ (the actual commencement time of the tipping) as well as the duration and character of the transient tipping trajectory and the character of the final equilibrium state. Such details are not available from extrapolation of a trend in stability. They require knowledge about global dynamics and may include surprises such as ``late tipping'' within a transient \cite{bastiaansen2023climate}. These dynamics may contain considerable uncertainties and surprises as other systems may be affected and models may be poorly parametrised in these areas. These details are vital in estimating the possible costs and benefits of any decisions made in response to a prediction.

Suppose that at time $t_p$ we predict tipping within the horizon $[t_p, t_p +\Delta]$, as in Section~\ref{subsec:extrap-trend}, and write $t_{\tip} \in [t_p, t_p +\Delta]$ for the time at which the threshold is first crossed. The \emph{lead time}  $\delta = t_{\tip} -t_p$ is the time available to a decision before a crossing. Any decision-making process will also introduce delays between issuing a warning and deciding on possible actions. This includes a time $\tau_{\mathrm{dec}}$ for a decision to be reached on how to act and a time $\tau_{\mathrm{act}}$ for the action to be implemented before it starts to take effect. Note that we count impact that may be due to anticipatory effects, for example in financial systems: knowledge of a future tipping can trigger impacts to be felt long before the tipping even happens.

\begin{figure}
    \centering
    \includegraphics[width=11cm]{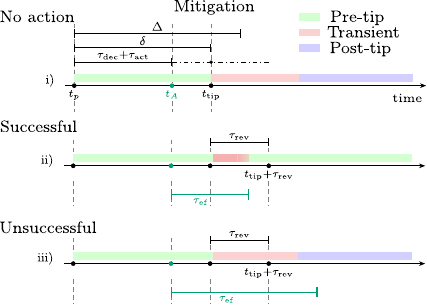}
    \caption{Schematic illustration of scenarios for a mitigation action to avoid predicted tipping. The system tips as time runs from left to right, from a light green ``Pre-tip" equilibrium to a light purple ``Post-tip" equilibrium, via a pink ``Transient". We show cases of i) no action, ii) successful and iii) unsuccessful mitigation action. A warning is given at time $t_p$ on a time horizon $\Delta$, and the actual start of tipping is at time $t_{tip}$. 
    There is a decision and action delay $\tau_{dec}+\tau_{act}$ before the time of action $t_A$.
    For ii) successful mitigation, the time of action $t_{ef}$ is such that the mitigation action is sufficiently strong within the period of duration $\tau_{rev}$ where it is possible to reverse commitment to tipping. For iii) unsuccessful mitigation, the action is too weak and/or comes too late to prevent tipping.}
    \label{fig:lead-times-mit}
\end{figure}

For a mitigation action (see Figure~\ref{fig:lead-times-mit}) the action needs to be early and strong enough to reverse the stability trend and keep/return the system to its pre-tip state. Even if the threshold is crossed, panel ii) shows that reversing the forcing quickly enough (shown as $\tau_{\mathrm{rev}}$ in the figure) can still avert the transition where a \emph{safe overshoot} is available \cite{ritchie2021overshooting}. On the other hand, if the action is not strong enough, the mitigation is unsuccessful as shown in panel iii), and the costs of the action may not bring the expected benefits.

\begin{figure}
    \centering
    \includegraphics[width=11cm]{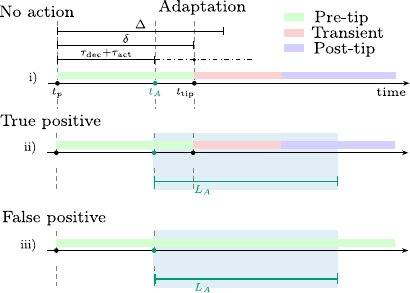}
    \caption{Schematic illustration of scenarios for an adaptation action, top panel and notation as in Figure~\ref{fig:lead-times-mit}. In this case, an adaptation action starts at time $t_A$ and runs for duration $L_A$. Typically, such an action will have initial costs but come with later benefits. We show cases of i) no action, ii) adaptation with a true positive prediction, and iii) adaptation with a false positive prediction. 
    For the true positive case ii), the system tips, but the adaptation action (light blue) reduces future impact. For the false positive case iii), the adaptation action is applied, but no tipping occurs.}
    \label{fig:lead-times-adapt}
\end{figure}

For an adaptation action (see Figure~\ref{fig:lead-times-adapt}) the action does not attempt to avoid the tipping but rather prepares for it. We show in light blue an adaptation action that will typically bring costs initially but benefits in future as the transient and post-tip scenario unfolds. Panel ii) shows the scenario of a true positive prediction, where the adaptation action starts before the transient. Panel iii) shows a scenario where the prediction is a false positive and the adaptation starts even though no tipping occurs.

These costs and benefits of possible actions should ideally be weighed against the benefits of avoiding the costs associated with the impacts of the tipping event. In the ideal ``no regrets'' case, where a mitigation or adaptation action carries no costs but only benefits, then an early warning system is not needed - the action should be implemented on insurance grounds.
In practice, even when mitigation is unsuccessful, it can potentially delay or extend the time at which one tips, and this can influence the time at which the impacts are felt.

\section{Discussion: improving skill and utility}
\label{sec:improving-skill}

In this section, we discuss approaches to improving skill, with reference to aspects highlighted by the colour shading in Figure~\ref{fig:skill-blobs}. We work though the potential for added skill though improved Earth System Models, improved use of observations, better decision-making processes and better integration procedures.

\subsection{Improving skill through better Earth System Models across the hierarchy}

Earth system models (ESMs), shown in pink in Figure~\ref{fig:skill-blobs}, are a vital component of our proposed forecasting system, especially for forecasting tipping events and regimes that are outside of historical observations and understanding the physical mechanisms behind these processes. 
As we highlighted previously in Section \ref{subsec:UQ}, some of the issues here can stem from errors in initial state estimation and parametric uncertainty, as well as structural errors. Hence improvements in the physics, parametrisation and numerical accuracy will be one route to improving skill of an early warning system.

The failure to acknowledge and account for structural errors in ESMs and other climate models can lead to biased estimates of parameters of interest and overconfident predictions \cite{brynjarsdottir_learning_2014}. Some of the future scenarios correspond to unprecedented states (i.e. those that have not been observed in the recorded history), which is ultimately an extrapolation from known physics and observations. There could be other physical processes not well represented in these models, or represented differentially across different models, such as the oceanic processes in the North Atlantic \cite{falkena2024causal}. 

Historical and observational data are crucial for addressing errors associated with model-based inference. For example, reconstruction of the past states is conditional to our assumptions about the physical processes at play at that time and invariably has uncertainty in the estimated time period and local environment. Some of these issues are being addressed by collecting data from more isolated and/or previously unexplored locations, advanced dating methods, and combining knowledge from multiple sources. These new observations can be used to guide future model development. We discuss this further in Sections \ref{subsec:existing_obs} and \ref{subsec:improve-skill-new-obs}.

However, there is only so much that can be done using historical/current observational data, particularly if the tipping phenomena of interest are not sampled in this data. This is where long-duration simulations, benchmarked on paleoclimate evidence of tipping phenomena, are particularly important \cite{Boers2022PaleoTipping}. For example, TipMIP \cite{tipmip2025} is creating a range of future forcing scenarios in ESMs where emissions are ramped up/maintained/ramped down and where various tipping phenomena can be studied.

Another way to navigate structural uncertainty is by testing the underlying physical mechanisms in low-order and medium complexity models. Critical outputs from these models should match with outputs from models of intermediate complexity and up in the hierarchy. Having low-order conceptual models (and emulators) can also help improve the catalogue of possible scenarios. Since complex models are expensive, one can at least test out hypotheses about future forcing to a crude approximation. 

To address parametric uncertainty in ESMs and other complex models, the multivariate ensembles produced by these models can be compared to the observed trends in data using approaches from inverse UQ \cite{kennedy2001bayesian, stuart2010inverse} to estimate parameters of interest and states. 
Moreover, they can help identify the appropriate observable(s) with high signal to noise ratio that represent current state of the system. This is related to {\bf H3} from section \ref{sec:math-challenges}.

\subsection{Improvements through better use of existing data}
\label{subsec:existing_obs}

As highlighted grey in Figure~\ref{fig:skill-blobs}, observations form the basis of the forecasting framework. Hence, making better use of the data available will be one route to improving the skill and utility of an early warning system, alongside the new observation streams discussed in Section~\ref{subsec:improve-skill-new-obs}. 
Getting more from these observations starts with knowing how much each stream contributes to skill. Observing System Experiments (OSE) \cite{magnussen2013,bormann2019} quantify this by comparing forecasts made with the full observing system with those made after withholding a particular stream. 

Observations and Earth system models are closely linked, so improvements in one can lead to improvements in the other. Comparing observation streams with one another and with model output reveals they disagree. Where independent records and modelled dynamics agree, confidence in the forecast is strengthened, while disagreement prompts further investigation, as the issue may lie with the observations, the model, or both \cite{simpson2025confronting}.
Records collected for other purposes may also carry information about stability that is missing or less clearly resolved in the data currently used. Recent use of bivalve data serves as an excellent example \cite{arellano-nava_historians_2025}.

Further gains come can come from how the records are analysed, whether by separating trends from artefacts introduced by changes in measurement noise or in how records from different instruments are merged \cite{Smith2023ReliabilityResilienceEstimation}, or by combining measured variables into an observable that shows critical slowing down more clearly than any individual record \cite{lohmann_choosing_2025}. Recent approaches also couple machine learning with data assimilation, where assimilation provides the state estimates used to train a machine-learning model, which then serves as the forecasting model in subsequent assimilation \cite{brajard_combining_2020}. Such approaches may allow early warnings to be developed from existing records even where the dynamics are not understood well enough to specify a low-order model in advance.

\subsection{Improvements through new observation streams}
\label{subsec:improve-skill-new-obs}

While it is clear that better use of existing observational data has the potential to improve forecasting skill, new observation streams have the potential to strengthen predictive skill through capturing crucial data that may transform an unpredictable phenomenon into one that is predictable. This includes additional observations that de-noise existing observations, and observations ones that hold entirely new signal information. This is especially the case for  ocean or cryosphere tipping elements where deep ocean or subsurface ice sheet measurements are not routinely collected, for obvious reasons, while weather balloons are regularly used to gain information that can be assimilated into circulation models.

What sort of additional observations are able to make decisive improvements to skill? This includes choice of environmental variables, as well as where and when to measure.
As additional observations may be technologically challenging, such questions can be considered using  Observing System Simulation Experiments (OSSE) \cite{kubo_predictability_2025}, where a potential new observing system is proposed to extract information from an Earth system model, and this can be compared to a ``ground truth" from the model. This includes trying to evaluate where predicted trends are observed \cite{simpson2025confronting}. Such experiments allow one to explore the effectiveness of active and adaptive observations.

It is particularly challenging to predict the timing of tipping events associated with rare events (noise-induced), for example, a glacial lake draining after an extreme rainfall event. Rare event algorithms can selectively produce high-impact, low-probability events from existing ESMs or intermediate complexity models, but it may be that the representative rare events will be unavailable for observation. In such cases, the best one can hope for is to evaluate the probabilities of the event occurring in a given time period.

\subsection{Improving skill through better communication and decision-making}

As highlighted in blue in Figure~\ref{fig:skill-blobs}, the utility of an early warning depends not only on the skill of the forecast, but also on the decision available to the relevant decision-maker and the consequences of acting or not acting \cite{leomoine2014watch}. This requires understanding the available portfolio of responses and how it may change over time. New technical developments may, for example, expand the range of mitigation and adaptation options, although they may introduce additional uncertainties of their own. The consequences of decisions must also be considered across true/false positives/negatives, and for different prediction horizons. These consequences may differ substantially between stakeholders and over different timescales, so their cost/benefit analysis will require considerable interdisciplinary effort.

An early warning system should therefore be developed with the needs of its intended users in mind \cite{popova2026designing}. Working with stakeholders to co-produce decision-relevant information falls within the broad area of climate services \cite{Lemos2005,Stocker2026}, where an important aim is to reduce the gap between scientific information and its practical usability; an often iterative process \cite{bremer2017}. The delays discussed in Section~\ref{subsec:tipping-duration-timing} must also be incorporated into the design of the warning system. In some settings, this may motivate a rapid communication route in which limited but decision-critical information is issued first, followed later by a more detailed assessment.

Such communication must make the uncertainties and limitations of the forecast clear. This relates to the argument of \cite{majszak2026metaphors} that the precautionary principle, commonly adopted in public discussion about tipping points, can undermine scientific credibility and informed decision-making around climate tipping points. 
Transparent communication of uncertainties in probabilities of events and their associated risks can guard against this  \cite{rietkerk2026warning}.  Guidance such as the AQuA Book \cite{analysisfunction2025aqua} similarly emphasises rigorous treatment of uncertainty and clear communication of the limitations of the analysis used to support the decisions.

Even with improved forecasts and communication, decisions concerning tipping points may have to be made under deep uncertainty. Approaches such as Adaptive Policymaking, Adaptive Pathways and Dynamic Adaptive Policy Pathways (DAPP) have been developed to support decisions that remain robust while retaining flexibility such that responses can be adjusted as conditions change \cite{walker2013adapt,haasnoot2013dynamic}. Early warning information can contribute to this process by updating knowledge of the current state, the approach towards critical conditions, and the uncertainty surrounding future evolution. Forecasting can also support related approaches such as backcasting, which identifies adaptation and mitigation actions by working backward from possible future states to the present \cite{biesbroek_adaptation_2025}.

\subsection{Improving skill through optimisation of Early Warning Systems}

Given models, observation data and decision-maker input, the yellow parts of Figure~\ref{fig:skill-blobs} highlight a number of aspects of an Early Warning System that are open to improvement. This includes better algorithms and methods of state and trend estimation, better specification of warning system design parameters (see \cite{huntingford2026}) and better ways to quantify, communicate and reduce uncertainty. For this, one of the most important yet challenging problems is to train EWS on suitable ensembles that include the full range of expected behaviours. There will be a trade-off between timescales and skill of prediction. For this, we suggest it will be important to understand potential impacts and to be led by these rather than specific Earth system events.

\subsection{Concluding remarks}

There are many open questions that are raised by the preceding discussion of ways to improve utility and skill of Early Warning Systems for tipping behaviour. Some of these are in areas of the physical sciences, mathematics and time series analysis, but many concern the social and economic impact of how decisions are made, who is making them, who is affected by them and over what timescale we wish to understand benefits and costs.

Moreover, the range of options available for action will change with time, either due to external factors (e.g. due to geopolitical events) or due to innovations in possible actions (e.g. due to innovations in technology or research). This means that the design parameters of the system will typically change with time. Moreover, different stakeholders will clearly have starkly different benefits. Hence, for effective decision-making in climate mitigation and adaptation, it will remain important for consensus and negotiation \cite{IPCC_AR6} to spread costs and benefits in an equitable way. 

For specific tipping elements in the Earth's climate system, one can often say more clearly which aspects in the schematic in Figure~\ref{fig:skill-blobs} are most important to prioritise for improvement. For example, tipping of the SPG (subpolar gyre) or AMOC involves regional or large-scale changes in ocean circulation where large parts of the system are poorly observed. In such cases, we expect that closing the deep-water observation gap will be most helpful to improve skill. For tipping elements such as boreal forest under threat of dieback, understanding the modelling of extreme atmospheric events and the modelling of dynamic vegetation may be the most helpful. Nonetheless, it is also clearly important to limit prediction to the predictable, though it is often not clear {\em a priori} whether this is possible!

\section*{Acknowledgements}

We thank all members of the AdvanTip team and from other ARIA projects ``Forecasting Tipping Points'' for stimulating discussions.
All authors thank the Advanced Research and Invention Agency (ARIA) for support via the AdvanTip project (grant no. SCOP-PR01-P003). PA, PR and TL thank the ClimTip project which has received funding from the European Union's Horizon Europe research and innovation programme under grant agreement No. 101137601: Funded by the European Union. PA thanks the 
Past to Future (P2F) project which has received funding from the European Union’s Horizon Europe research and innovation programme under grant agreement No. 101184070: Funded by the European Union. Views and opinions expressed are however those of the author(s) only and do not necessarily reflect those of the European Union or the European Climate, Infrastructure and Environment Executive Agency (CINEA). Neither the European Union nor the granting authority can be held responsible for them. 
For the purpose of open access, the authors have applied a Creative Commons Attribution (CC BY) licence to any Author Accepted Manuscript version arising from this submission.

\bibliographystyle{plain}
\bibliography{ews-refs-cleaned}

@article{hurlstone2025threshold,
    author = {Hurlstone, Mark J. and White, Ben and Newell, Ben R.},
    title = {Threshold uncertainty, early warning signals and the prevention of dangerous climate change},
    journal = {Royal Society Open Science},
    volume = {12},
    number = {3},
    pages = {240425},
    year = {2025},
    month = {03},
    issn = {2054-5703},
    doi = {10.1098/rsos.240425},
    url = {https://doi.org/10.1098/rsos.240425},
    eprint = {https://royalsocietypublishing.org/rsos/article-pdf/doi/10.1098/rsos.240425/18531/rsos.240425.pdf},
}

@article{barrett2022sensitivity,
    author={Barrett, S. and Dannenberg, A.},
    title={Sensitivity of collective action to uncertainty about climate tipping points},
    journal={Nature Clim Change},
    volume={4}, 
    pages={36–39},
    year={2014},
    doi={https://doi.org/10.1038/nclimate2059}
}

@article{leomoine2014watch,
Author = {Lemoine, Derek and Traeger, Christian},
Title = {Watch Your Step: Optimal Policy in a Tipping Climate},
Journal = {American Economic Journal: Economic Policy},
Volume = {6},
Number = {1},
Year = {2014},
Month = {February},
Pages = {137–66},
DOI = {10.1257/pol.6.1.137}
}

@article{ferro2012bias,
author = {Ferro, C. A. T. and Fricker, T. E.},
title = {A bias-corrected decomposition of the Brier score},
journal = {Quarterly Journal of the Royal Meteorological Society},
volume = {138},
number = {668},
pages = {1954-1960},
doi = {https://doi.org/10.1002/qj.1924},
url = {https://rmets.onlinelibrary.wiley.com/doi/abs/10.1002/qj.1924},
eprint = {https://rmets.onlinelibrary.wiley.com/doi/pdf/10.1002/qj.1924},
year = {2012}
}

@article{huntingford2026,
author={Chris Huntingford and al},
title={Ideas for Creating Early Warning Systems in the Climate Sciences},
journal={preprint},
year={2026}
}

@article{tipsci2026,
author={Gabi Hegerl and al},
title={High-impact climate events and tipping dynamics: assessing scientific understanding and future challenges},
journal={preprint},
year={2026}
}

@misc{bormann2019,
  author = {Niels Bormann and Heather Lawrence and Jacqueline Farnan},
  title = {Global observing system experiments in the ECMWF assimilation system},
  year = {2019},
  journal = {ECMWF Technical Memoranda},
  number = {839},
  month = {01/2019},
  publisher = {ECMWF},
  url = {https://www.ecmwf.int/node/18859},
  doi = {10.21957/sr184iyz},
  language = {eng},
}

@article {magnussen2013,
      author = "Linus Magnusson and Erland Källén",
      title = "Factors Influencing Skill Improvements in the ECMWF Forecasting System",
      journal = "Monthly Weather Review",
      year = "2013",
      publisher = "American Meteorological Society",
      address = "Boston MA, USA",
      volume = "141",
      number = "9",
      doi = "10.1175/MWR-D-12-00318.1",
      pages=      "3142 - 3153",
      url = "https://journals.ametsoc.org/view/journals/mwre/141/9/mwr-d-12-00318.1.xml"
}

@Article{tipmip2025,
AUTHOR = {Winkelmann, R. and Dennis, D. P. and Donges, J. F. and Loriani, S. and Klose, A. K. and Abrams, J. F. and Alvarez-Solas, J. and Albrecht, T. and Armstrong McKay, D. and Bathiany, S. and Blasco Navarro, J. and Brovkin, V. and Burke, E. and Danabasoglu, G. and Donner, R. V. and Dr\"uke, M. and Georgievski, G. and Goelzer, H. and Harper, A. B. and Hegerl, G. and Hirota, M. and Hu, A. and Jackson, L. C. and Jones, C. and Kim, H. and Koenigk, T. and Lawrence, P. and Lenton, T. M. and Liddy, H. and Lic\'on-Sal\'aiz, J. and Menthon, M. and Montoya, M. and Nitzbon, J. and Nowicki, S. and Otto-Bliesner, B. and Pausata, F. and Rahmstorf, S. and Ramin, K. and Robinson, A. and Rockstr\"om, J. and Romanou, A. and Sakschewski, B. and Sch\"adel, C. and Sherwood, S. and Smith, R. S. and Steinert, N. J. and Swingedouw, D. and Willeit, M. and Weijer, W. and Wood, R. and Wyser, K. and Yang, S.},
TITLE = {The Tipping Points Modelling Intercomparison Project (TIPMIP): Assessing tipping point risks in the Earth system},
JOURNAL = {EGUsphere},
VOLUME = {2025},
YEAR = {2025},
PAGES = {1--52},
URL = {https://egusphere.copernicus.org/preprints/2025/egusphere-2025-1899/},
DOI = {10.5194/egusphere-2025-1899}
}

@article{stuart2010inverse,
  title={Inverse problems: a Bayesian perspective},
  author={Stuart, Andrew M},
  journal={Acta numerica},
  volume={19},
  pages={451--559},
  year={2010},
  publisher={Cambridge University Press}
}

@book{owen2001empirical,
  title={Empirical likelihood},
  author={Owen, Art B and others},
  volume={92},
  year={2001},
  publisher={Chapman \& Hall/CRC New York}
}

@article{stefanski2002calculus,
  title={The calculus of M-estimation},
  author={Stefanski, Leonard A and Boos, Dennis D},
  journal={The American Statistician},
  volume={56},
  number={1},
  pages={29--38},
  year={2002},
  publisher={Taylor \& Francis}
}

@article{Lemos2005,
  author    = {Lemos, Maria Carmen and Morehouse, Barbara J.},
  title     = {The co-production of science and policy in integrated climate assessments},
  journal   = {Global Environmental Change},
  volume    = {15},
  number    = {1},
  pages     = {57--68},
  year      = {2005},
  publisher = {Elsevier},
  doi       = {10.1016/j.gloenvcha.2004.09.004},
  url       = {https://doi.org/10.1016/j.gloenvcha.2004.09.004}
}

@article{bremer2017,
author = {Bremer, Scott and Meisch, Simon},
title = {Co-production in climate change research: reviewing different perspectives},
journal = {WIREs Climate Change},
volume = {8},
number = {6},
pages = {e482},
doi = {https://doi.org/10.1002/wcc.482},
url = {https://wires.onlinelibrary.wiley.com/doi/abs/10.1002/wcc.482},
eprint = {https://wires.onlinelibrary.wiley.com/doi/pdf/10.1002/wcc.482},
year = {2017}
}

@article{Stocker2026,
  author  = {Stocker, Thomas F. and Belcher, Stephen E. and Brunet, Gilbert and
             Dube, Pauline and Fr{\"o}licher, Thomas L. and Hegerl, Gabriele C. and
             Jain, Shipra and Marengo, Jose and Metternicht, Graciela and
             Naish, Tim and van der Wel, Narelle and Vautard, Robert},
  title   = {Insights from Climate Tipping Point Science for the Forthcoming
             {IPCC} Seventh Assessment Report and Reflections for National
             Meteorological and Hydrological Services},
  journal = {Bulletin of the American Meteorological Society},
  year    = {2026},
  month   = jul,
  volume  = {107},
  number  = {7},
  pages   = {E1579--E1585},
  doi     = {10.1175/BAMS-D-26-0119.1}
}

@article{bissiri2016general,
  title={A general framework for updating belief distributions},
  author={Bissiri, Pier Giovanni and Holmes, Chris C and Walker, Stephen G},
  journal={Journal of the Royal Statistical Society Series B: Statistical Methodology},
  volume={78},
  number={5},
  pages={1103--1130},
  year={2016},
  publisher={Oxford University Press}
}

@article{matsubara2022robust,
  title={Robust generalised Bayesian inference for intractable likelihoods},
  author={Matsubara, Takuo and Knoblauch, Jeremias and Briol, Fran{\c{c}}ois-Xavier and Oates, Chris J},
  journal={Journal of the Royal Statistical Society Series B: Statistical Methodology},
  volume={84},
  number={3},
  pages={997--1022},
  year={2022},
  publisher={Oxford University Press}
}

@article{kriegler2009imprecise,
  title={Imprecise probability assessment of tipping points in the climate system},
  author={Kriegler, Elmar and Hall, Jim W and Held, Hermann and Dawson, Richard and Schellnhuber, Hans Joachim},
  journal={Proceedings of the national Academy of Sciences},
  volume={106},
  number={13},
  pages={5041--5046},
  year={2009},
  publisher={National Academy of Sciences}
}

@article{lam2022climate,
  title={Climate tipping points and expert judgment},
  author={Lam, Vincent and Majszak, Mason M},
  journal={Wiley Interdisciplinary Reviews: Climate Change},
  volume={13},
  number={6},
  pages={e805},
  year={2022},
  publisher={Wiley Online Library}
}

@article{majszak2026metaphors,
  title={Metaphors and values in the communication of climate tipping points},
  author={Majszak, Mason and Martini, Carlo},
  journal={Synthese},
  volume={207},
  number={3},
  pages={101},
  year={2026},
  publisher={Springer}
}

@article{haasnoot2013dynamic,
  title={Dynamic adaptive policy pathways: A method for crafting robust decisions for a deeply uncertain world},
  author={Haasnoot, Marjolijn and Kwakkel, Jan H and Walker, Warren E and Ter Maat, Judith},
  journal={Global environmental change},
  volume={23},
  number={2},
  pages={485--498},
  year={2013},
  publisher={Elsevier}
}

@article{kimpton2025challenges,
  title={Challenges and opportunities in uncertainty quantification for healthcare and biological systems},
  author={Kimpton, Louise M and Paun, L Mihaela and Colebank, Mitchel J and Volodina, Victoria},
  journal={Philosophical transactions. Series A, Mathematical, physical, and engineering sciences},
  volume={383},
  number={2292},
  pages={20240232},
  year={2025}
}

@book{cox2006principles,
  title={Principles of statistical inference},
  author={Cox, David Roxbee},
  year={2006},
  publisher={Cambridge university press}
}

@misc{analysisfunction2025aqua,
  author       = {{Government Analysis Function}},
  title        = {The {AQuA} {Book}},
  howpublished = {GOV.UK guidance},
  year         = {2025},
  month        = jul,
  note         = {Revised edition, published 30 July 2025; last updated 20 July 2026.
                  Produced with the Government Operational Research Service, Government
                  Actuary's Department and Office for National Statistics},
  url          = {https://www.gov.uk/guidance/the-aqua-book},
  urldate      = {2026-09-11}
}

@article{goldstein2013assessing,
  title={Assessing model adequacy},
  author={Goldstein, Michael and Seheult, Allan and Vernon, Ian},
  journal={Environmental modelling: finding simplicity in complexity},
  pages={435--449},
  year={2013},
  publisher={Wiley Online Library}
}

@article{goldstein2025systematic,
  title={Systematic structural discrepancy assessment for computer models},
  author={Goldstein, Michael and Vernon, Ian and Cumming, Jonathan A},
  journal={Philosophical transactions. Series A, Mathematical, physical, and engineering sciences},
  volume={383},
  number={2293},
  pages={20240214},
  year={2025}
}

@article{jacques2024estimation,
  title={Estimation of AMOC transition probabilities using a machine learning--based rare-event algorithm},
  author={Jacques-Dumas, Val{\'e}rian and van Westen, Ren{\'e} M and Dijkstra, Henk A},
  journal={Artificial Intelligence for the Earth Systems},
  volume={3},
  number={4},
  pages={e240002},
  year={2024},
  publisher={American Meteorological Society}
}

@article{castellana2019transition,
  title={Transition probabilities of noise-induced transitions of the Atlantic Ocean circulation},
  author={Castellana, Daniele and Baars, Sven and Wubs, Fred W and Dijkstra, Henk A},
  journal={Scientific Reports},
  volume={9},
  number={1},
  pages={20284},
  year={2019},
  publisher={Nature Publishing Group UK London}
}

@article{bellomo2021future,
  title={Future climate change shaped by inter-model differences in Atlantic meridional overturning circulation response},
  author={Bellomo, Katinka and Angeloni, Michela and Corti, Susanna and von Hardenberg, Jost},
  journal={Nature Communications},
  volume={12},
  number={1},
  pages={3659},
  year={2021},
  publisher={Nature Publishing Group UK London}
}

@article{mehling2024limits,
  title={Limits to predictability of the asymptotic state of the Atlantic Meridional Overturning Circulation in a conceptual climate model},
  author={Mehling, Oliver and B{\"o}rner, Reyk and Lucarini, Valerio},
  journal={Physica D: Nonlinear Phenomena},
  volume={459},
  pages={134043},
  year={2024},
  publisher={Elsevier}
}

@article{cini2024simulating,
  title={Simulating AMOC tipping driven by internal climate variability with a rare event algorithm},
  author={Cini, Matteo and Zappa, Giuseppe and Ragone, Francesco and Corti, Susanna},
  journal={npj Climate and Atmospheric Science},
  volume={7},
  number={1},
  pages={31},
  year={2024},
  publisher={Nature Publishing Group UK London}
}

@article{van2023asymmetry,
  title={Asymmetry of AMOC hysteresis in a state-of-the-art global climate model},
  author={van Westen, Ren{\'e} M and Dijkstra, Henk A},
  journal={Geophysical Research Letters},
  volume={50},
  number={22},
  pages={e2023GL106088},
  year={2023},
  publisher={Wiley Online Library}
}

@article{caiado2015bayesian,
  title={Bayesian uncertainty analysis for complex physical systems modelled by computer simulators with applications to tipping points},
  author={Caiado, Camila CS and Goldstein, M},
  journal={Communications in Nonlinear Science and Numerical Simulation},
  volume={26},
  number={1-3},
  pages={123--136},
  year={2015},
  publisher={Elsevier}
}

@article{cini2026assessing,
  title={Assessing the probability of CO$_2$ - driven AMOC transitions using a rare event algorithm in PlaSim-LSG},
  author={Cini, Matteo and Dijkstra, Henk A and Jacques-Dumas, Val{\'e}rian and Zappa, Giuseppe},
  journal={Journal of Physics: Complexity},
  volume={7},
  number={2},
  pages={02LT02},
  year={2026},
  publisher={IOP Publishing}
}

@article{volodina2021importance,
  title={The importance of uncertainty quantification in model reproducibility},
  author={Volodina, Victoria and Challenor, Peter},
  journal={Philosophical transactions. Series A, Mathematical, physical, and engineering sciences},
  volume={379},
  number={2197},
  pages={20200071},
  year={2021}
}

@article{kennedy2001bayesian,
  title={Bayesian calibration of computer models},
  author={Kennedy, Marc C and O'Hagan, Anthony},
  journal={Journal of the Royal Statistical Society: Series B (Statistical Methodology)},
  volume={63},
  number={3},
  pages={425--464},
  year={2001},
  publisher={Wiley Online Library}
}

@article{lucarini2024detecting,
  title={Detecting and attributing change in climate and complex systems: {F}oundations, {G}reen’s functions, and nonlinear fingerprints},
  author={Lucarini, Valerio and Chekroun, Micka{\"e}l D},
  journal={Physical Review Letters},
  volume={133},
  number={24},
  pages={244201},
  year={2024},
  publisher={APS}
}

@article{chekroun2020ruelle,
  title={Ruelle--Pollicott Resonances of Stochastic Systems in Reduced State Space. Part I: Theory: MD Chekroun et al.},
  author={Chekroun, Micka{\"e}l D and Tantet, Alexis and Dijkstra, Henk A and Neelin, J David},
  journal={Journal of Statistical Physics},
  volume={179},
  number={5},
  pages={1366--1402},
  year={2020},
  publisher={Springer}
}

@article{zagli2026bridging,
  title={Bridging the gap between Koopmanism and response theory: Using natural variability to predict forced response},
  author={Zagli, Niccol{\`o} and Colbrook, Matthew J and Lucarini, Valerio and Mezi{\'c}, Igor and Moroney, John},
  journal={SIAM Journal on Applied Dynamical Systems},
  volume={25},
  number={1},
  pages={196--229},
  year={2026},
  publisher={SIAM}
}

@book{kutz2016dynamic,
  title={Dynamic mode decomposition: data-driven modeling of complex systems},
  author={Kutz, J Nathan and Brunton, Steven L and Brunton, Bingni W and Proctor, Joshua L},
  year={2016},
  publisher={SIAM}
}

@article{williams2015data,
  title={A data--driven approximation of the koopman operator: Extending dynamic mode decomposition},
  author={Williams, Matthew O and Kevrekidis, Ioannis G and Rowley, Clarence W},
  journal={Journal of Nonlinear Science},
  volume={25},
  number={6},
  pages={1307--1346},
  year={2015},
  publisher={Springer}
}

@article{walker2013adapt,
  title={Adapt or perish: A review of planning approaches for adaptation under deep uncertainty},
  author={Walker, Warren E and Haasnoot, Marjolijn and Kwakkel, Jan H},
  journal={Sustainability},
  volume={5},
  number={3},
  pages={955--979},
  year={2013},
  publisher={MDPI}
}

@article{schmid2022dynamic,
  title={Dynamic mode decomposition and its variants},
  author={Schmid, Peter J},
  journal={Annual Review of Fluid Mechanics},
  volume={54},
  number={1},
  pages={225--254},
  year={2022},
  publisher={Annual Reviews}
}

@article{kwasniok2024semiparametric,
  author={Frank Kwasniok},
  title={Semiparametric maximum likelihood reconstruction of stochastic differential equations driven by white and correlated noise},
  journal={Chaos: An Interdisciplinary Journal of Nonlinear Science},
  volume={34},
  number={10},
  year={2024},
  publisher={AIP Publishing}
}

@article{rietkerk2026warning,
  title={Warning for climate tipping without compromising credibility},
  author={Rietkerk, Max and Haasnoot, Marjolijn and Vuuren, Detlef P van},
  journal={Climatic Change},
  volume={179},
  number={6},
  pages={122},
  year={2026},
  publisher={Springer}
}

@article{boettiger2012early,
  title={Early warning signals and the prosecutor's fallacy},
  author={Boettiger, Carl and Hastings, Alan},
  journal={Proceedings of the Royal Society B: Biological Sciences},
  volume={279},
  number={1748},
  pages={4734},
  year={2012}
}

@article{rietkerk2025ambiguity,
  title={Ambiguity of early warning signals for climate tipping points},
  author={Rietkerk, Max and Skiba, Vanessa and Weinans, Els and H{\'e}bert, Rapha{\"e}l and Laepple, Thomas},
  journal={Nature Climate Change},
  volume={15},
  number={5},
  pages={479--488},
  year={2025},
  publisher={Nature Publishing Group UK London}
}

@article{simpson2025confronting,
author = {Isla R. Simpson  and Tiffany A. Shaw  and Paulo Ceppi  and Amy C. Clement  and Erich Fischer  and Kevin M. Grise  and Angeline G. Pendergrass  and James A. Screen  and Robert C. J. Wills  and Tim Woollings  and Russell Blackport  and Joonsuk M. Kang  and Stephen Po-Chedley },
title = {Confronting Earth System Model trends with observations},
journal = {Science Advances},
volume = {11},
number = {11},
pages = {eadt8035},
year = {2025},
doi = {10.1126/sciadv.adt8035},
URL = {https://www.science.org/doi/abs/10.1126/sciadv.adt8035},
eprint = {https://www.science.org/doi/pdf/10.1126/sciadv.adt8035}
}

@book{lutkepohl2013introduction,
  title={Introduction to multiple time series analysis},
  author={L{\"u}tkepohl, Helmut},
  year={2013},
  publisher={Springer Science \& Business Media}
}

@book{gardiner2004handbook,
  title={Handbook of stochastic methods},
  author={Gardiner, Crispin W.},
  volume={3},
  year={2004},
  publisher={springer Berlin}
}

@article{pilipovic2025strang,
  title={Strang splitting for parametric inference in second-order stochastic differential equations},
  author={Pilipovic, Predrag and Samson, Adeline and Ditlevsen, Susanne},
  journal={Stochastic Processes and their Applications},
  volume={187},
  pages={104650},
  year={2025},
  publisher={Elsevier}
}

@article{falkena2024causal,
  title={Causal Mechanisms of Subpolar Gyre Variability in CMIP6 Models},
  author={Falkena, Swinda KJ and Dijkstra, Henk A and von der Heydt, Anna S},
  journal={arXiv preprint arXiv:2408.16541},
  year={2024}
}

@article{ritchie2021overshooting,
  title={Overshooting tipping point thresholds in a changing climate},
  author={Ritchie, Paul DL and Clarke, Joseph J and Cox, Peter M and Huntingford, Chris},
  journal={Nature},
  volume={592},
  number={7855},
  pages={517--523},
  year={2021},
  publisher={Nature Publishing Group}
}

@article{ashwin2025early,
  title={Early warning skill, extrapolation and tipping for accelerating cascades},
  author={Ashwin, Peter and Bastiaansen, Robbin and von der Heydt, Anna S and Ritchie, Paul},
  journal= {Proc. Roy. Soc. Lond. A},
  volume={To appear},
  pages={},
  url={https://arxiv.org/abs/2506.01981},
  year={2025}
}

@article{wieczorek2023rate,
  title={Rate-induced tipping: Thresholds, edge states and connecting orbits},
  author={Wieczorek, Sebastian and Xie, Chun and Ashwin, Peter},
  journal={Nonlinearity},
  volume={36},
  number={6},
  pages={3238},
  year={2023},
  publisher={IOP Publishing}
}

@article{bastiaansen2023climate,
  title={Climate response and sensitivity: time scales and late tipping points},
  author={Bastiaansen, Robbin and Ashwin, Peter and von der Heydt, Anna S},
  journal={Proceedings of the Royal Society A},
  volume={479},
  number={2269},
  pages={20220483},
  year={2023},
  publisher={The Royal Society}
}

@article{Boers2022PaleoTipping,
	title = {Theoretical and paleoclimatic evidence for abrupt transitions in the {Earth} system},
	volume = {17},
	doi = {10.1088/1748-9326/ac8944},
	number = {9},
	journal = {Environmental Research Letters},
	author = {Boers, Niklas and Ghil, Michael and Stocker, Thomas F},
	month = sep,
	year = {2022},
}

@article{deb2022machine,
  title={Machine learning methods trained on simple models can predict critical transitions in complex natural systems},
  author={Deb, Smita and Sidheekh, Sahil and Clements, Christopher F and Krishnan, Narayanan C and Dutta, Partha S},
  journal={Royal Society Open Science},
  volume={9},
  number={2},
  pages={211475},
  year={2022}
}

@inproceedings{batchnorm,
  title={Batch normalization: Accelerating deep network training by reducing internal covariate shift},
  author={Ioffe, Sergey and Szegedy, Christian},
  booktitle={International conference on machine learning},
  pages={448--456},
  year={2015},
  organization={pmlr}
}

@article{santurkar2018does,
  title={How does batch normalization help optimization?},
  author={Santurkar, Shibani and Tsipras, Dimitris and Ilyas, Andrew and Madry, Aleksander},
  journal={Advances in neural information processing systems},
  volume={31},
  year={2018}
}

@article{sevinchan2026context,
  title={In-context learning to predict critical transitions in dynamical systems},
  author={Sevinchan, Yunus and Nathaniel, Juan and Ueltzh{\"o}ffer, Kai and Roesch, Carla and Weber, Tobias and Laschos, Vaios and Fan, Hang and Ramien, Gregor and Haux, Johannes and Gentine, Pierre and others},
  journal={arXiv preprint arXiv:2605.12308},
  year={2026}
}

@article{zhang2025anticipating,
  title={Anticipating {AMOC} transitions via deep learning},
  author={Zhang, Wenjie and Huang, Yu and Bathiany, Sebastian and Shin, Yechul and Ben-Yami, Maya and Zhou, Suiping and Boers, Niklas},
  journal={arXiv preprint arXiv:2509.06450},
  year={2025}
}

@article{guardamagna2026predicting,
  title={Predicting the distance of the {AMOC} to its tipping point using CNNs},
  author={Guardamagna, Francesco and Sinet, Sacha and Dijkstra, Henk A},
  journal={EGUsphere},
  volume={2026},
  pages={1--48},
  year={2026},
  publisher={Copernicus Publications G{\"o}ttingen, Germany}
}

@article{Boettner2022CSDCorrNoise,
	title = {Critical slowing down in dynamical systems driven by nonstationary correlated noise},
	volume = {4},
	doi = {10.1103/PhysRevResearch.4.013230},
	number = {1},
	journal = {Physical Review Research},
	author = {Boettner, Christopher and Boers, Niklas},
	month = mar,
	year = {2022},
}

@article{BenYami2023AMOCDataCSD,
	title = {Uncertainties in critical slowing down indicators of observation-based fingerprints of the {Atlantic} {Overturning} {Circulation}},
	volume = {14},
	issn = {2041-1723},
	doi = {10.1038/s41467-023-44046-9},
	number = {1},
	journal = {Nature Communications},
	author = {Ben-Yami, Maya and Skiba, Vanessa and Bathiany, Sebastian and Boers, Niklas},
	month = dec,
	year = {2023},
	pages = {8344},
}

@article{Smith2023ReliabilityResilienceEstimation,
	title = {Reliability of resilience estimation based on multi-instrument time series},
	volume = {14},
	issn = {2190-4987},
	doi = {10.5194/esd-14-173-2023},
	number = {1},
	journal = {Earth System Dynamics},
	author = {Smith, Taylor and Zotta, Ruxandra-Maria and Boulton, Chris A. and Lenton, Timothy M. and Dorigo, Wouter and Boers, Niklas},
	month = feb,
	year = {2023},
	pages = {173--183},
}

@article{Proverbio2022COVIDEWS,
	title = {Performance of early warning signals for disease re-emergence: {A} case study on {COVID}-19 data},
	volume = {18},
	issn = {1553-7358},
	shorttitle = {Performance of early warning signals for disease re-emergence},
	doi = {10.1371/journal.pcbi.1009958},
	number = {3},
	journal = {PLOS Computational Biology},
	author = {Proverbio, Daniele and Kemp, Françoise and Magni, Stefano and Gonçalves, Jorge},
	editor = {Althouse, Benjamin},
	month = mar,
	year = {2022},
	pages = {e1009958},
}

@book{Berglund2006NoiseInducedPhenomena,
	series = {Probability and {Its} {Applications}},
	title = {Noise-{Induced} {Phenomena} in {Slow}-{Fast} {Dynamical} {Systems}},
	isbn = {978-1-84628-038-2},
	publisher = {Springer-Verlag},
	author = {Berglund, Nils},
	year = {2006},
	doi = {10.1007/1-84628-186-5},
}

@article{laitinen2021probabilistic,
  title={Probabilistic early warning signals},
  author={Laitinen, Ville and Dakos, Vasilis and Lahti, Leo},
  journal={Ecology and Evolution},
  volume={11},
  number={20},
  pages={14101--14114},
  year={2021},
  publisher={Wiley Online Library}
}

@article{ives2012detecting,
  title={Detecting dynamical changes in nonlinear time series using locally linear state-space models},
  author={Ives, Anthony R and Dakos, Vasilis},
  journal={Ecosphere},
  volume={3},
  number={6},
  pages={1--15},
  year={2012},
  publisher={Wiley Online Library}
}

@inproceedings{laitinen2022probabilistic,
  title={Probabilistic multivariate early warning signals},
  author={Laitinen, Ville and Lahti, Leo},
  booktitle={International Conference on Computational Methods in Systems Biology},
  pages={259--274},
  year={2022},
  organization={Springer}
}

@article{Morr2024InternalNoiseInterference,
author = {Morr, Andreas and Boers, Niklas and Ashwin, Peter},
title = {Internal Noise Interference to Warnings of Tipping Points in Generic Multidimensional Dynamical Systems},
journal = {SIAM Journal on Applied Dynamical Systems},
volume = {23},
number = {4},
pages = {2793-2806},
year = {2024},
doi = {10.1137/24M1669104},
URL = {https://doi.org/10.1137/24M1669104}
}

@book{Kuznetsov2004AppliedBifurcationTheory,
	series = {Applied {Mathematical} {Sciences}},
	title = {Elements of {Applied} {Bifurcation} {Theory}},
	volume = {112},
	isbn = {978-1-4419-1951-9 978-1-4757-3978-7},
	publisher = {Springer New York},
	author = {Kuznetsov, Yuri A.},
	year = {2004},
	doi = {10.1007/978-1-4757-3978-7},
}

@article{Williamson2015MultiDBifDetection,
	title = {Detection of bifurcations in noisy coupled systems from multiple time series},
	volume = {25},
	issn = {1054-1500, 1089-7682},
	doi = {10.1063/1.4908603},
	number = {3},
	journal = {Chaos: An Interdisciplinary Journal of Nonlinear Science},
	author = {Williamson, Mark S. and Lenton, Timothy M.},
	month = mar,
	year = {2015},
	pages = {036407},
}

@article{Ditlevsen2023AMOCPrediction,
	title = {Warning of a forthcoming collapse of the {Atlantic} meridional overturning circulation},
	volume = {14},
	doi = {10.1038/s41467-023-39810-w},
	number = {1},
	journal = {Nature Communications},
	author = {Ditlevsen, Peter and Ditlevsen, Susanne},
	month = jul,
	year = {2023},
	pages = {4254},
}

@article{gneiting2007strictly,
  title={Strictly proper scoring rules, prediction, and estimation},
  author={Gneiting, Tilmann and Raftery, Adrian E},
  journal={Journal of the American statistical Association},
  volume={102},
  number={477},
  pages={359--378},
  year={2007},
  publisher={Taylor \& Francis}
}

@article{Dai2015StabilityvsResilience,
	title = {Relation between stability and resilience determines the performance of early warning signals under different environmental drivers},
	volume = {112},
	doi = {10.1073/pnas.1418415112},
	number = {32},
	journal = {Proceedings of the National Academy of Sciences},
	author = {Dai, Lei and Korolev, Kirill S. and Gore, Jeff},
	year = {2015},
	pages = {10056--10061},
}

@article{Hastings2010NoEWS,
	title = {Regime shifts in ecological systems can occur with no warning},
	volume = {13},
	doi = {10.1111/j.1461-0248.2010.01439.x},
	number = {4},
	journal = {Ecology Letters},
	author = {Hastings, Alan and Wysham, Derin B.},
	year = {2010},
	pages = {464--472},
}

@article{Chorin2000MZMemory,
	title = {Optimal prediction and the {Mori}–{Zwanzig} representation of irreversible processes},
	volume = {97},
	doi = {10.1073/pnas.97.7.2968},
	number = {7},
	journal = {Proceedings of the National Academy of Sciences},
	author = {Chorin, Alexandre J. and Hald, Ole H. and Kupferman, Raz},
	year = {2000},
	pages = {2968--2973},
}

@article{Zwanzig1961Memory,
	title = {Memory effects in irreversible thermodynamics},
	volume = {124},
	doi = {10.1103/PhysRev.124.983},
	number = {4},
	journal = {Physical Review},
	author = {Zwanzig, Robert},
	month = nov,
	year = {1961},
	pages = {983--992},
}

@article{Brett2020ROCEWS,
	title = {Detecting critical slowing down in high-dimensional epidemiological systems},
	volume = {16},
	doi = {10.1371/journal.pcbi.1007679},
	number = {3},
	journal = {PLOS Computational Biology},
	author = {Brett, Tobias and Ajelli, Marco and Liu, Quan-Hui and Krauland, Mary G. and Grefenstette, John J. and van Panhuis, Willem G. and Vespignani, Alessandro and Drake, John M. and Rohani, Pejman},
	month = mar,
	year = {2020},
	pages = {1--19},
}

@article{Thompson2010FoldTippingPredict,
	title = {Climate tipping as a noisy bifurcation: a predictive technique},
	volume = {76},
	issn = {0272-4960},
	doi = {10.1093/imamat/hxq060},
	number = {1},
	journal = {IMA Journal of Applied Mathematics},
	author = {Thompson, J. M. T. and Sieber, Jan},
	month = dec,
	year = {2010},
	pages = {27--46},
}

@article{Boerlijst2013SilentCollapse,
	title = {Catastrophic collapse can occur without early warning: {Examples} of silent catastrophes in structured ecological models},
	volume = {8},
	doi = {10.1371/journal.pone.0062033},
	number = {4},
	journal = {PLOS ONE},
	author = {Boerlijst, Maarten C. and Oudman, Thomas and de Roos, André M.},
	month = apr,
	year = {2013},
	pages = {1--6},
}

@article{Kuehn2011CSD,
	title = {A mathematical framework for critical transitions: {Bifurcations}, fast–slow systems and stochastic dynamics},
	volume = {240},
	issn = {0167-2789},
	doi = {10.1016/j.physd.2011.02.012},
	number = {12},
	journal = {Physica D: Nonlinear Phenomena},
	author = {Kuehn, Christian},
	year = {2011},
	pages = {1020--1035},
}

@article{Kuehn2013CSDVar,
	title = {A mathematical framework for critical transitions: {Normal} forms, variance and applications},
	volume = {23},
	doi = {10.1007/s00332-012-9158-x},
	number = {3},
	journal = {Journal of Nonlinear Science},
	author = {Kuehn, Christian},
	month = jun,
	year = {2013},
	pages = {457--510},
}

@article{krakovska2024resilience,
  title={Resilience of dynamical systems},
  author={Krakovsk{\'a}, Hana and Kuehn, Christian and Longo, Iacopo P},
  journal={European Journal of Applied Mathematics},
  volume={35},
  number={1},
  pages={155--200},
  year={2024},
  publisher={Cambridge University Press}
}

@article{morr2026computing,
  title={Computing resilience measures in dynamical systems},
  author={Morr, Andreas and Kuehn, Christian and Datseris, George},
  journal={Chaos: An Interdisciplinary Journal of Nonlinear Science},
  volume={36},
  number={2},
  year={2026},
  publisher={AIP Publishing}
}

@article{smith2026estimating,
  title={Estimating the Resilience of Non-Stationary Systems},
  author={Smith, Taylor and Morr, Andreas and Sch{\"o}tz, Christof and Boers, Niklas},
  journal={arXiv preprint arXiv:2604.24345},
  year={2026}
}

@article{lehnertz2024time,
  title={Time-series-analysis-based detection of critical transitions in real-world non-autonomous systems},
  author={Lehnertz, Klaus},
  journal={Chaos: An Interdisciplinary Journal of Nonlinear Science},
  volume={34},
  number={7},
  year={2024},
  publisher={AIP Publishing}
}

@article{
richardson2023,
author = {Katherine Richardson  and Will Steffen  and Wolfgang Lucht  and Jørgen Bendtsen  and Sarah E. Cornell  and Jonathan F. Donges  and Markus Drüke  and Ingo Fetzer  and Govindasamy Bala  and Werner von Bloh  and Georg Feulner  and Stephanie Fiedler  and Dieter Gerten  and Tom Gleeson  and Matthias Hofmann  and Willem Huiskamp  and Matti Kummu  and Chinchu Mohan  and David Nogués-Bravo  and Stefan Petri  and Miina Porkka  and Stefan Rahmstorf  and Sibyll Schaphoff  and Kirsten Thonicke  and Arne Tobian  and Vili Virkki  and Lan Wang-Erlandsson  and Lisa Weber  and Johan Rockström },
title = {Earth beyond six of nine planetary boundaries},
journal = {Science Advances},
volume = {9},
number = {37},
pages = {eadh2458},
year = {2023},
doi = {10.1126/sciadv.adh2458},
URL = {https://www.science.org/doi/abs/10.1126/sciadv.adh2458},
eprint = {https://www.science.org/doi/pdf/10.1126/sciadv.adh2458},
}

@article{fawcett2006introduction,
  title={An introduction to {ROC} analysis},
  author={Fawcett, Tom},
  journal={Pattern recognition letters},
  volume={27},
  number={8},
  pages={861--874},
  year={2006},
  publisher={Elsevier}
}

@article{dakos2008slowing,
  title={Slowing down as an early warning signal for abrupt climate change},
  author={Dakos, Vasilis and Scheffer, Marten and Van Nes, Egbert H. and Brovkin, Victor and Petoukhov, Vladimir and Held, Hermann},
  journal={Proceedings of the National Academy of Sciences},
  volume={105},
  number={38},
  pages={14308--14312},
  year={2008},
  publisher={National Academy of Sciences}
}

@article{scheffer2012anticipating,
  title={Anticipating critical transitions},
  author={Scheffer, Marten and Carpenter, Stephen R. and Lenton, Timothy M. and Bascompte, Jordi and Brock, William and Dakos, Vasilis and Van de Koppel, Johan and Van de Leemput, Ingrid A. and Levin, Simon A. and Van Nes, Egbert H. and others},
  journal={science},
  volume={338},
  number={6105},
  pages={344--348},
  year={2012},
  publisher={American Association for the Advancement of Science}
}

@article{dablander2022deep,
  title={Deep learning for tipping points: preprocessing matters},
  author={Dablander, Fabian and Bury, Thomas M},
  journal={Proceedings of the National Academy of Sciences},
  volume={119},
  number={37},
  pages={e2207720119},
  year={2022},
  publisher={National Academy of Sciences}
}

@article{bury2021deep,
  title={Deep learning for early warning signals of tipping points},
  author={Bury, Thomas M. and Sujith, R.I. and Pavithran, Induja and Scheffer, Marten and Lenton, Timothy M. and Anand, Madhur and Bauch, Chris T.},
  journal={Proceedings of the National Academy of Sciences},
  volume={118},
  number={39},
  pages={e2106140118},
  year={2021},
  publisher={National Academy of Sciences}
}

@article{dakos2024tipping,
  title={Tipping point detection and early warnings in climate, ecological, and human systems},
  author={Dakos, Vasilis and Boulton, Chris A. and Buxton, Joshua E and Abrams, Jesse F. and Arellano-Nava, Beatriz and Armstrong McKay, David I. and Bathiany, Sebastian and Blaschke, Lana and Boers, Niklas and Dylewsky, Daniel and others},
  journal={Earth System Dynamics},
  volume={15},
  number={4},
  pages={1117--1135},
  year={2024},
  publisher={Copernicus GmbH}
}

@article{armstrong2022exceeding,
author = {David I. Armstrong McKay  and Arie Staal  and Jesse F. Abrams  and Ricarda Winkelmann  and Boris Sakschewski  and Sina Loriani  and Ingo Fetzer  and Sarah E. Cornell  and Johan Rockström  and Timothy M. Lenton },
title = {Exceeding 1.5°C global warming could trigger multiple climate tipping points},
journal = {Science},
volume = {377},
number = {6611},
pages = {eabn7950},
year = {2022},
doi = {10.1126/science.abn7950},
URL = {https://www.science.org/doi/abs/10.1126/science.abn7950}}

@article{ashwin2017parameter,
  title={Parameter shifts for nonautonomous systems in low dimension: bifurcation-and rate-induced tipping},
  author={Ashwin, Peter and Perryman, Clare and Wieczorek, Sebastian},
  journal={Nonlinearity},
  volume={30},
  number={6},
  pages={2185},
  year={2017},
  publisher={IOP Publishing}
}

@article{chapman2025quantifying,
    author = {Chapman, Ruth R. and Ashwin, Peter and Wood, Richard A.},
    title = {Quantifying tipping behavior: Geometric early warnings and quasipotentials for a box model of AMOC},
    journal = {Chaos: An Interdisciplinary Journal of Nonlinear Science},
    volume = {35},
    number = {2},
    pages = {023164},
    year = {2025},
    month = {02},
    issn = {1054-1500},
    doi = {10.1063/5.0240959},
    url = {https://doi.org/10.1063/5.0240959},
    eprint = {https://pubs.aip.org/aip/cha/article-pdf/doi/10.1063/5.0240959/20413754/023164_1_5.0240959.pdf},
}

@article{kuehn2013mathematical,
  title={A mathematical framework for critical transitions: normal forms, variance and applications},
  author={Kuehn, Christian},
  journal={Journal of Nonlinear Science},
  volume={23},
  pages={457--510},
  year={2013},
  publisher={Springer},
doi = {10.1007/s00332-012-9158-x}
}

@article{
benmayi2024uncertainties,
author = {Maya Ben-Yami  and Andreas Morr  and Sebastian Bathiany  and Niklas Boers},
title = {Uncertainties too large to predict tipping times of major Earth system components from historical data},
journal = {Science Advances},
volume = {10},
number = {31},
pages = {eadl4841},
year = {2024},
doi = {10.1126/sciadv.adl4841},
URL = {https://www.science.org/doi/abs/10.1126/sciadv.adl4841},
}

@article{dakos2014critical,
  title={Critical slowing down as early warning for the onset of collapse in mutualistic communities},
  author={Dakos, Vasilis and Bascompte, Jordi},
  journal={Proceedings of the National Academy of Sciences},
  volume={111},
  number={49},
  pages={17546--17551},
  year={2014},
doi = {10.1073/pnas.1406326111}
}

@article{ashwin2012tipping,
  title={Tipping points in open systems: bifurcation, noise-induced and rate-dependent examples in the climate system},
  author={Ashwin, Peter and Wieczorek, Sebastian and Vitolo, Renato and Cox, Peter},
  journal={Philosophical Transactions of the Royal Society A: Mathematical, Physical and Engineering Sciences},
  volume={370},
  number={1962},
  pages={1166--1184},
  year={2012},
  publisher={The Royal Society Publishing},
doi = {10.1098/rsta.2011.0306}
}

@article{stommel61,
author = {Stommel,H.},
title = {Thermohaline Convection with Two Stable Regimes of Flow},
journal = {Tellus},
volume = {13},
number = {2},
pages = {224-230},
year  = {1961},
publisher = {Taylor & Francis},
doi = {10.1111/j.2153-3490.1961.tb00079.x}
}

@article{huang2024anticipating,
  title={Deep Learning for predicting rate-induced tipping},
  author={Huang, Yu and Bathiany, Sebastian and Ashwin, Peter and Boers, Niklas},
  year={2024},
  journal={Nature Machine Intelligence},
  volume={6}, 
  pages={1556–1565}
}

@article{ashwin2025contrasting,
author = {Ashwin, Peter and Newman, Julian and R\"{o}mer, Raphael},
title = {Contrasting Chaotic and Stochastic Forcing: Tipping Windows and Attractor Crises},
journal = {SIAM Journal on Applied Dynamical Systems},
volume = {24},
number = {1},
pages = {277-316},
year = {2025},
doi = {10.1137/24M1661534},
URL = {https://doi.org/10.1137/24M1661534},
}

@article{Lenton08,
author={Lenton, Timothy M. and Held, Hermann and Kriegler, Elmar and Hall, Jim W. and Lucht, Wolfgang and Rahmstorf, Stefan and Schellnhuber, Hans Joachim},
title = {Tipping elements in the {Earth}'s climate system},
journal = {Proceedings of the National Academy of Sciences},
volume = {105},
number = {6},
pages = {1786-1793},
year = {2008},
doi = {10.1073/pnas.0705414105},
}

@article{brunton2016discovering,
  title={Discovering governing equations from data by sparse identification of nonlinear dynamical systems},
  author={Brunton, Steven L and Proctor, Joshua L and Kutz, J Nathan},
  journal={Proceedings of the national academy of sciences},
  volume={113},
  number={15},
  pages={3932--3937},
  year={2016},
  publisher={National Academy of Sciences}
}

@article{kevrekidis2009equation,
  title={Equation-free multiscale computation: Algorithms and applications},
  author={Kevrekidis, Ioannis G and Samaey, Giovanni},
  journal={Annual review of physical chemistry},
  volume={60},
  number={1},
  pages={321--344},
  year={2009},
  publisher={Annual Reviews}
}

@article{morr2026extrapolation,
  title={Extrapolation from historical data cannot reliably predict the time of a potential AMOC collapse},
  author={Morr, Andreas and Ben-Yami, Maya and Groenke, Brian and Sch{\"o}tz, Christof and Cotronei, Alessandro and Myrvoll-Nilsen, Eirik and Bathiany, Sebastian and Rypdal, Martin and Boers, Niklas},
  journal={arXiv preprint arXiv:2604.20341},
  year={2026}
}

@misc{ritchie2026evaluating,
      title={Evaluating the skill of a geometric early warning for tipping in a rapidly forced nonlinear system}, 
      author={Paul D. L. Ritchie and Sneha Kachhara and Peter Ashwin},
      year={2026},
      eprint={2605.16128},
      archivePrefix={arXiv},
      primaryClass={math.DS},
      url={https://arxiv.org/abs/2605.16128}, 
}

@book{IPCC_AR6,
author = {IPCC and Masson-Delmotte, V. and Zhai, P. and Pirani, A. and Connors, S.L. and Péan, C. and Berger, S. and Caud, N. and Chen, Y. and Goldfarb, L. and Gomis, M.I. and Huang, M. and Leitzell, K. and Lonnoy, E. and Matthews, J.B.R. and Maycock, T.K. and Waterfield, T. and Yelekçi, O. and Yu, R. and Zhou, B.},
title = {{IPCC, 2021: Climate Change 2021: The Physical Science Basis. Contribution of Working Group I to the Sixth Assessment Report of the Intergovernmental Panel on
Climate Change}},
publisher = {Cambridge University Press, Cambridge, United Kingdom and New York, NY, USA, },
doi = {doi:10.1017/9781009157896.},
year = {2021},
}

@article{lohmann2025role,
  title={The role of edge states for early warning of tipping points},
  author={Lohmann, Johannes and Hansen, Alfred B. and Lovo, Alessandro and Chapman, Ruth and Bouchet, Freddy and Lucarini, Valerio},
  journal={Proceedings of the Royal Society A: Mathematical, Physical and Engineering Sciences},
  volume={481},
  number={2312},
  year={2025},
  publisher={The Royal Society}
}

@article{hobden2026early,
author={Bryony Hobden and Paul D. Ritchie and Peter Ashwin},
title={Early warnings of critical transitions through
vector autoregression: lessons from multiscale
systems},
journal={arXiv},
url={https://doi.org/10.48550/arXiv.2605.28260},
year={2026}
}

@article{melbourne2011note,
   title={A note on diffusion limits of chaotic skew-product flows},
   volume={24},
   ISSN={1361-6544},
   url={http://dx.doi.org/10.1088/0951-7715/24/4/018},
   DOI={10.1088/0951-7715/24/4/018},
   number={4},
   journal={Nonlinearity},
   publisher={IOP Publishing},
   author={Melbourne, I and Stuart, A M},
   year={2011},
   month=mar, pages={1361–1367}
}

@book{gottwald2017stochastic,
author={Georg A. Gottwald and Daan Crommelin and Christian Franzke},
title={Stochastic climate theory},
booktitle={Nonlinear and Stochastic Climate Dynamics}, 
publisher={Cambridge University Press},
year={2017}
}

@article{Hasselmann1993optimal,
      author = "K.  Hasselmann",
      title = "Optimal Fingerprints for the Detection of Time-dependent Climate Change",
      journal = "Journal of Climate",
      year = "1993",
      publisher = "American Meteorological Society",
      address = "Boston MA, USA",
      volume = "6",
      number = "10",
      doi = "10.1175/1520-0442(1993)006<1957:OFFTDO>2.0.CO;2",
      pages=      "1957 - 1971",
      url = "https://journals.ametsoc.org/view/journals/clim/6/10/1520-0442_1993_006_1957_offtdo_2_0_co_2.xml"
}

@article{zhang2015predictability,
  title = {Predictability of critical transitions},
  author = {Zhang, Xiaozhu and Kuehn, Christian and Hallerberg, Sarah},
  journal = {Phys. Rev. E},
  volume = {92},
  issue = {5},
  pages = {052905},
  numpages = {14},
  year = {2015},
  month = {Nov},
  publisher = {American Physical Society},
  doi = {10.1103/PhysRevE.92.052905}
}

@article{kelly2017deterministic,
title = {Deterministic homogenization for fast–slow systems with chaotic noise},
journal = {J. Funct. Anal.},
volume = {272},
number = {10},
pages = {4063--4102},
year = {2017},
author = {David~Kelly and Ian~Melbourne}
}

@article{diekert2025doearly,
    author = {Diekert, Florian and Heyen, Daniel and Nesje, Frikk and Shayegh, Soheil},
    title = {Do early warning signals of tipping points lead to better decisions?},
    journal = {Journal of The Royal Society Interface},
    volume = {22},
    number = {225},
    pages = {20240864},
    year = {2025},
    month = {04},
    issn = {1742-5689},
    doi = {10.1098/rsif.2024.0864},
    url = {https://doi.org/10.1098/rsif.2024.0864},
    eprint = {https://royalsocietypublishing.org/rsif/article-pdf/doi/10.1098/rsif.2024.0864/2830314/rsif.2024.0864.pdf},
}

@article{myrvoll-nilsen_bayesian_2025,
    title = {Bayesian analysis of early warning signals using a time-dependent model},
    volume = {16},
    issn = {2190-4979},
    url = {https://esd.copernicus.org/articles/16/1539/2025/},
    doi = {10.5194/esd-16-1539-2025},
    language = {en-US},
    number = {5},
    urldate = {2026-05-07},
    journal = {Earth System Dynamics},
    publisher = {Copernicus GmbH},
    author = {Myrvoll-Nilsen, Eirik and Hallali, Luc and Rypdal, Martin},
    month = sep,
    year = {2025},
    pages = {1539--1556},
}

@article{biesbroek_adaptation_2025,
    title = {Adaptation planning in the context of a weakening and possibly collapsing {Atlantic} {Meridional} {Overturning} {Circulation} ({AMOC})},
    volume = {25},
    issn = {1436-378X},
    url = {https://doi.org/10.1007/s10113-025-02434-5},
    doi = {10.1007/s10113-025-02434-5},
    language = {en},
    number = {3},
    urldate = {2026-05-23},
    journal = {Regional Environmental Change},
    author = {Biesbroek, Robbert and Haasnoot, Marjolijn and Mach, Katharine J. and Petersen, Arthur C.},
    month = jul,
    year = {2025},
    pages = {93},
}

@article{kuhlmann_seizure_2018,
    title = {Seizure prediction — ready for a new era},
    volume = {14},
    copyright = {2018 Springer Nature Limited},
    issn = {1759-4766},
    url = {https://www.nature.com/articles/s41582-018-0055-2},
    doi = {10.1038/s41582-018-0055-2},
    language = {en},
    number = {10},
    urldate = {2026-06-08},
    journal = {Nature Reviews Neurology},
    publisher = {Nature Publishing Group},
    author = {Kuhlmann, Levin and Lehnertz, Klaus and Richardson, Mark P. and Schelter, Björn and Zaveri, Hitten P.},
    month = oct,
    year = {2018},
    pages = {618--630},
}

@article{machete_early_2013,
    title = {Early {Warning} with {Calibrated} and {Sharper} {Probabilistic} {Forecasts}},
    volume = {32},
    issn = {1099-131X},
    url = {https://onlinelibrary.wiley.com/doi/abs/10.1002/for.2242},
    doi = {10.1002/for.2242},
    language = {en},
    number = {5},
    urldate = {2026-06-10},
    journal = {Journal of Forecasting},
    author = {Machete, Reason L.},
    year = {2013},
    note = {\_eprint: https://onlinelibrary.wiley.com/doi/pdf/10.1002/for.2242},
    pages = {452--468},
}

@techreport{woolliams_general_2025,
    address = {Teddington, Middlesex},
    type = {Technical {Report}},
    title = {General guidance on a metrological approach to fundamental data records ({FDR}), thematic data products ({TDP}) and fiducial reference measurements ({FRM}) – {Uncertainty} {Analysis} {Process}},
    url = {https://qa4eo.org/docs/3_Process_Document.pdf},
    institution = {National Physical Laboratory},
    author = {Woolliams, Emma and Behnia, Sajedeh and Bialek, Agniezka and De Vis, Pieter and Fahy, Jacob and Fox, Nigel and Gardiner, Tom and Green, Paul and Hunt, Samuel and Mittaz, Jonathan and Mota, Bernardo and Nightingale, Joane and Origo, Niall and Pustogvar, Anna and Stedman, Madeline},
    year = {2025},
}

@article{popova2026designing,
  title={Designing decision-centric early warning systems for climate tipping points: the North Atlantic Subpolar Gyre as a test case},
  author={Popova, Ekaterina and Bai, Xinyi and Clark, Elisabeth and Cooper, Charlotte and Gommenginger, Christine and Jacobs, Zoe and Jude, Simon and Kelly, Stephen and Loveridge, Alexandra and Marzocchi, Alice and others},
  year={2026},
  publisher={EarthArXiv},
  journal={Preprint}
}

@book{gelman_bayesian_2013,
    address = {New York},
    edition = {3},
    title = {Bayesian {Data} {Analysis}},
    isbn = {978-0-429-11307-9},
    doi = {10.1201/b16018},
    publisher = {Chapman and Hall/CRC},
    author = {Gelman, Andrew and Carlin, John B. and Stern, Hal S. and Dunson, David B. and Vehtari, Aki and Rubin, Donald B.},
    month = nov,
    year = {2013},
}

@article{brynjarsdottir_learning_2014,
doi = {10.1088/0266-5611/30/11/114007},
url = {https://doi.org/10.1088/0266-5611/30/11/114007},
year = {2014},
month = {oct},
publisher = {IOP Publishing},
volume = {30},
number = {11},
pages = {114007},
author = {Brynjarsdottir, Jenny and O'Hagan, Anthony},
title = {Learning about physical parameters: the importance of model discrepancy},
journal = {Inverse Problems}
}

@article{boulton_early_2026,
    title = {Early warning signals of tipping points: {A} deep learning approach to direction and timing},
    copyright = {CC BY Attribution 4.0 International},
    shorttitle = {Early warning signals of tipping points},
    url = {https://eartharxiv.org/repository/view/14610/},
    language = {en},
    journal={Preprint},
    urldate = {2026-09-15},
    publisher = {EarthArXiv},
    author = {Boulton, Chris A. and Fadera, Muhammed},
    month = aug,
    year = {2026},
}

@article{kubo_predictability_2025,
    title = {Predictability of {Climate} {Tipping} {Focusing} on {Internal} {Variability} in the {Earth} {System}},
    volume = {52},
    issn = {1944-8007},
    url = {https://onlinelibrary.wiley.com/doi/abs/10.1029/2024GL113146},
    doi = {10.1029/2024GL113146},
    language = {en},
    number = {11},
    urldate = {2026-09-21},
    journal = {Geophysical Research Letters},
    author = {Kubo, Amane and Sawada, Yohei},
    year = {2025},
    pages = {e2024GL113146},
}

@article{arellano-nava_historians_2025,
    title = {From historians to forecasters: {The} potential of bivalve records to assess resilience and provide early warnings for marine tipping points},
    volume = {10},
    issn = {2378-2242},
    shorttitle = {From historians to forecasters},
    url = {https://onlinelibrary.wiley.com/doi/abs/10.1002/lol2.10455},
    doi = {10.1002/lol2.10455},
    language = {en},
    number = {2},
    urldate = {2026-09-20},
    journal = {Limnology and Oceanography Letters},
    author = {Arellano-Nava, Beatriz and Boulton, Chris A. and Reynolds, David J. and Butler, Paul G. and Scourse, James and Lenton, Timothy M. and Halloran, Paul R.},
    year = {2025},
    pages = {200--211},
}

@article{lohmann_choosing_2025,
    title = {Choosing observables that capture critical slowing down before tipping points: {A} {Fokker}-{Planck} operator approach},
    volume = {112},
    issn = {2470-0045, 2470-0053},
    shorttitle = {Choosing observables that capture critical slowing down before tipping points},
    url = {https://link.aps.org/doi/10.1103/l2v2-xndy},
    doi = {10.1103/l2v2-xndy},
    language = {en},
    number = {6},
    urldate = {2026-09-25},
    journal = {Physical Review E},
    author = {Lohmann, Johannes and Gottwald, Georg A.},
    month = dec,
    year = {2025},
    pages = {064204},
}

@article{brajard_combining_2020,
    title = {Combining data assimilation and machine learning to emulate a dynamical model from sparse and noisy observations: {A} case study with the {Lorenz} 96 model},
    volume = {44},
    issn = {1877-7503},
    shorttitle = {Combining data assimilation and machine learning to emulate a dynamical model from sparse and noisy observations},
    url = {https://www.sciencedirect.com/science/article/pii/S1877750320304725},
    doi = {10.1016/j.jocs.2020.101171},
    urldate = {2026-09-23},
    journal = {Journal of Computational Science},
    author = {Brajard, Julien and Carrassi, Alberto and Bocquet, Marc and Bertino, Laurent},
    month = jul,
    year = {2020},
    pages = {101171},
}

\end{document}